\documentclass[journal,twoside,web,a4paper]{IEEEtran}

\usepackage{amsthm}
\usepackage{generic}
\usepackage{cite}
\usepackage{amsmath,amssymb,amsfonts}
\usepackage{xcolor}
\usepackage{url}
\usepackage{algorithmic}
\usepackage{graphicx}
\usepackage{algorithm,algorithmic}
\usepackage{enumitem}
\usepackage{textcomp}
\usepackage{caption}
\usepackage{subcaption}
\newtheorem{assumption}{Assumption}
\newtheorem{corollary}{Corollary}
\newtheorem{theorem}{Theorem}

\newtheorem{lemma}{Lemma}
\newtheorem{proposition}{Proposition}
\newtheorem{example}{Example}
\newtheorem{definition}{Definition}
\newtheorem{case}{Case}
\allowdisplaybreaks
\def\BibTeX{{\rm B\kern-.05em{\sc i\kern-.025em b}\kern-.08em
    T\kern-.1667em\lower.7ex\hbox{E}\kern-.125emX}}
\begin{document}
\title{Convex Co-Design of Control Barrier Functions and Safe Feedback Controllers Under Input Constraints}
\author{Han Wang, Kostas Margellos, Antonis Papachristodoulou, and Claudio De Persis
\thanks{For the purpose of Open Access, the authors have applied a CC BY public copyright licence to any Author Accepted Manuscript (AAM) version arising from this submission. Part of the results of this manuscript has been submitted to the IEEE Conference on Decision and Control 2024. Here we significantly extend the conference version by additionally considering input constraints, adding new analysis examples, and generalizing the setting of safety specifications.}
\thanks{H. Wang, K. Margellos, and A. Papachristodoulou are with Department of Engineering Science, University of Oxford, OX1 3PJ, Oxford, UK. Emails: \{han.wang, kostas.margellos, antonis\}@eng.ox.ac.uk. }
\thanks{C. De Persis is with Engineering and Technology Institute, University of Groningen, 9747AG, The Netherlands. Email: c.de.persis@rug.nl.}
\thanks{AP was supported in part by UK's Engineering, Physical Sciences Research Council projects EP/X017982/1 and EP/Y014073/1.}}

\maketitle

\begin{abstract}
We study the problem of co-designing control barrier functions (CBF) and linear state feedback controllers for continuous-time linear systems. We achieve this by means of a single semi-definite optimization program. Our formulation can handle mixed-relative degree problems without requiring an explicit safe controller. Different $\mathcal{L}$-norm based input limitations can be introduced as convex constraints in the proposed program. We demonstrate our results on an omni-directional car numerical example.

\end{abstract}


\section{Introduction}
\label{sec:introduction}
\IEEEPARstart{S}{afety} is essential for feedback control systems. As a system is steered from an initial set to a target set, safety requires that the trajectory of the system avoids entering an unexpected region, or to remain inside a safe set. On the state space, safety is always formulated by means of constraints imposed on states. Based on these descriptions, two questions are raised: given a dynamical system $\dot x =f(x,u)$, a set of initial sets $\mathcal{I}$, and a set of safe states $\mathcal{S}$, (i) verify whether there exists a control input $u(\cdot)$, so that the trajectories starting from $\mathcal{I}$ stay inside $\mathcal{S}$; (ii) design such a control law $u(\cdot)$ that guarantees safety. The Control Barrier Functions (CBF) approach answers these two questions by using a continuously differentiable function that satisfies certain properties \cite{prajna2004safety,wieland2007constructive,ames2014control}.

A CBF aims to separate the safe and unsafe regions by its zero super- and sub-level sets; the initial set also belongs to the level set. In addition, there exists a control law, such that the vector field points towards the safe side on its zero sub-level set \cite{prajna2004safety}. This property is also known as \emph{invariance}, characterized by Nagumo's theorem \cite{nagumo1942lage}. It is therefore guaranteed that if the system starts from a point inside the zero super-level set, the system can always stay inside. Given a CBF, the controller that guarantees safety can be designed according to the direction requirement of vector field. However, synthesizing a CBF is not a trivial task even for linear systems. In general, even verifying a CBF is an NP-hard problem \cite[Proposition 2]{clark2022semi}. 

{
Designing a CBF is even more challenging when the relative degree between the function defines safe set and the system dynamics is high or mixed \cite{clark2021verification}. For relative degree we mean the number of times we need to differentiate a function whose level set encodes the safe set along the system dynamics until the control explicitly shows \cite{xiao2019control,tan2021high}. High or mixed relative degree is commonly seen in robotics collision avoidance problems, where the safe set is usually defined over positions for the obstacles, but the control signals are imposed on accelerations. 
} 

\subsection{Motivating Cases}\label{sec:case}
\begin{case}[Pathological vector field of CBF-QP]\label{cas:1}
Consider a continuous-time linear system with state matrix $A=\begin{bmatrix}
    -1&-1\\0&-1
\end{bmatrix}$ and input matrix $B=\begin{bmatrix}
    1\\1
\end{bmatrix}$. The system is controllable. Let $x=[x_1,x_2]^\top\in\mathbb{R}^2$ denote the states. The unsafe region is defined by $\mathcal{S}^c=\{x|s(x)\le 0\}$, where $s(x)=x_1^2+x_2^2-1$. Using $s(x)$ as a control barrier function in a quadratic programming framework \cite{ames2014control} for controller design, we obtain:
\begin{equation*}
    \begin{split}
        \min~&u^\top u,\\\mathrm{subject~to}~&\dot s(x)+10s(x)\ge 0.
    \end{split}
\end{equation*}
Following \cite{krstic2023inverse} , the analytical solution is given by 
\begin{equation*}
    u_s(x)=-\frac{\min\{0,f(x)\}}{g(x)},
\end{equation*}
where 
\begin{equation*}
    \begin{split}
        f(x)&=\frac{\partial s(x)}{\partial x}Ax+10s(x)\\
        &=8x_1^2-2x_1x_2+8x_2^2-10,\\
        g(x)&=\frac{\partial s(x)}{\partial x}B=2x_1+2x_2.
    \end{split}
\end{equation*}
\begin{figure*}[htbp]
     \centering
     \begin{subfigure}[t]{0.47\textwidth}
         \centering
         \includegraphics[width=\textwidth]{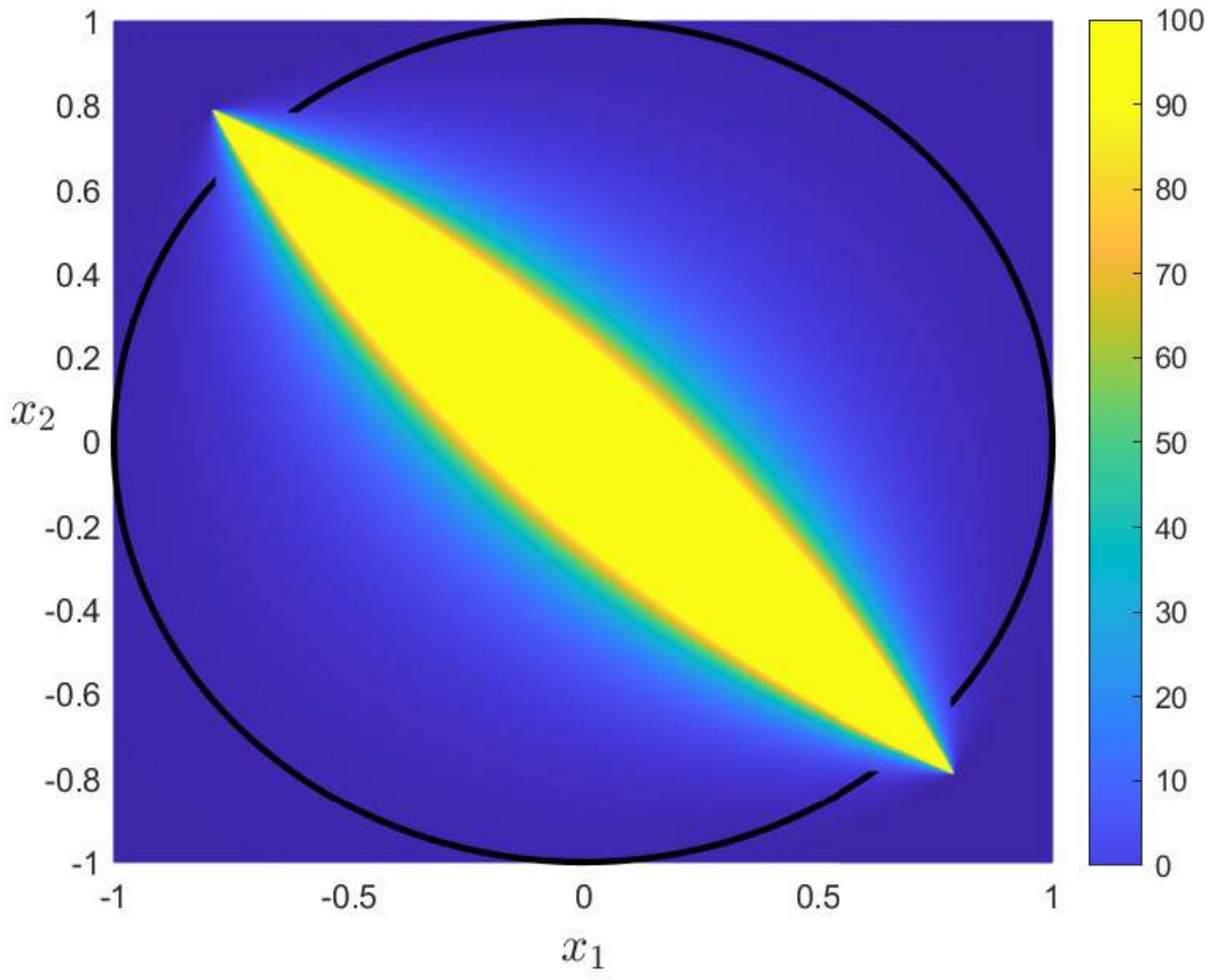}
         \caption{Values of $||u_s(x)||_2^2$ for $-1\le x_1 \le 1$, $-1\le x_2\le 1$. The value of $||u_s(x)||_2^2$ is limited to $100$ for visualization. The controller is only locally smooth, and the Lipschitz constant is large in a local region as the value varies a lot with little state changes. Different selection of a class-$\mathcal{K}$ function does not change the result for $x\in\partial \mathcal{B}$, which is the black curve in the figure. Clearly our designed feedback controller $u_b(x)=1.4164x_1 + 0.59702x_2$ is globally smooth as it is linear.}
         \label{fig:case1b}
     \end{subfigure}
     \hfill
     \begin{subfigure}[t]{0.47\textwidth}
         \centering
         \includegraphics[width=\textwidth]{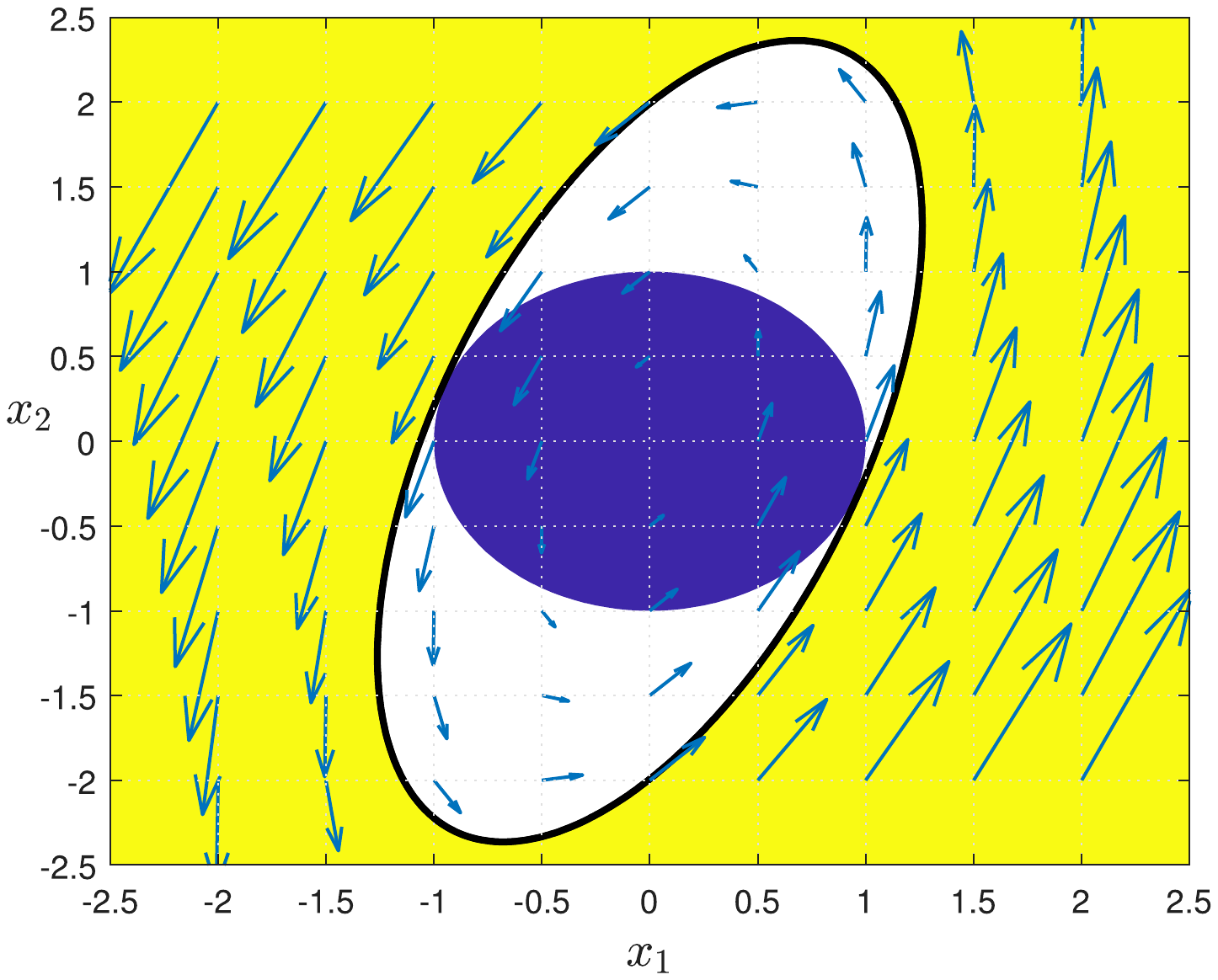}
         \caption{Comparison of the two control barrier functions in Case \ref{cas:1}. The blue round region is the unsafe set $\mathcal{S}^c$. The yellow open region is the control invariant set $\mathcal{B}:=\{x|b(x)\ge 0\}$, and $\partial \mathcal{B}:=\{x|b(x)= 0\}$ is the black curve. Blue arrows represents the vector field $Ax+Bu_b(x)$, which points inward $\mathcal{B}$ on $\partial \mathcal{B}$.}
         \label{fig:case1a}
     \end{subfigure}
        \caption{Simulation results for Case \ref{cas:1}.}
        \label{fig:case1}
\end{figure*}

When $g(x)$ tends to zero, and $f(x)<0$, $u_s(x)$ tends to infinity. As a consequence, the system cannot be safe at some points, especially points on $\partial \mathcal{S}$ with limited control authority. We also show in Figure \ref{fig:case1b} that the Lipschitz constant of $u_s(x)$ is very large. Later on, we will show that, by solving the proposed convex program \eqref{eq:dqcbflem} with $||u_b(x)||_2^2\le 8$ for $x$ such that $b(x)=0$, we obtain a new control barrier function $b(x)=0.88391x_1^2 - 0.50767x_1x_2 + 0.25205x_2^2 - 1$, and a feedback controller $u_b(x)= 1.4164x_1 + 0.59702x_2$. Comparison of the two control barrier functions is shown in Figure \ref{fig:case1a}.
It can be seen that the value of $||u_s(x)||_2^2$ is comparably large for small $x_2$. Meanwhile, our synthesized controller is constrained by $||u_b(x)||_2^2\le 8$ for $x\in\partial \mathcal{B}$.
\end{case}

\begin{case}[Mixed relative degree]\label{cas:2}
Consider a third-order continuous-time linear system with $\dot {x}_1=x_2+x_3$, $\dot {x}_2=x_1+u_1$, $\dot {x}_3=x_1+u_2$, where $x=[x_1,x_2,x_3]^\top \in\mathbb{R}^3$ is the state, $u=[u_1,u_2]^\top \in\mathbb{R}^2$ is the input. The unsafe region is defined by $\mathcal{S}^c:=\{x|s(x)\le 0\}$, where $s(x)=x_1^2+x_2^2-1$. { Let the relative degree be the number of times we need to
differentiate $s(x)$ along the dynamics until the control input $u$
appears in the resulting expression.
}For this case, the relative degree between $s(x)$ and the system is mixed, as the input $u_1$ appears in the first derivative of $s(x)$, whereas $u_2$ appears in the second derivative. $s(x)$ can not be directly used as a CBF using high-relative degree (exponential) CBF techniques \cite{xiao2019control,tan2021high,nguyen2016exponential}. 
By solving the convex program \eqref{eq:dqcbflem} that we will propose in the sequel, we obtain a control barrier function $b(x)=x_1^2+x_2^2-0.0129x_3^2-1$, and a feedback controller $u_1(x)=-2x_1+38.9x_2$, $u_2(x)=76.8x_1-0.5x_3$, which guarantees safety for the system. Clearly, the relative degree between $b(x)$ and the system dynamics is one. We highlight here that the backstepping CBF method \cite{taylor2022safe} would require a series of explicit pre-synthesized safe controllers which are, however, not needed for our method, which only requires the solution of a convex program.
    
\end{case}
\subsection{Related Work}
Dating back to the 1980's, there has been tremendous work on control invariance, especially for linear systems \cite{bertsekas1972infinite,kolmanovsky1998theory,castelan1993invariant,blanchini1999set,li2014complete,lin2007set}. For continuous-time linear systems, a half plane divided by an eigenvector is invariant \cite{castelan1993invariant}. For discrete-time systems, an invariant set can be constructed iteratively by state propagation. These methods focus on invariance but not safety. Building upon invariance, different methodologies have been proposed to synthesize control barrier functions. 

The first type of methods is reachability-based methods. Given a target set and a safe set, solving an optimal control problem returns a set of states starting from which the dynamical system can stay in the safe set and reach the target set. Such a set is usually the zero super-level set of a value function. Naturally, if only safety is considered over a finite horizon in the optimal control problem, the value function is a finite-time CBF \cite{choi2021robust}. More recently, the relationship between the safe value function and a CBF has been established \cite{massiani2022safe}. Solving this problem directly involves computing the solution of a Hamilton-Jacobi partial differential equation \cite{margellos2011hamilton,lygeros2004reachability,mitchell2005time,tomlin2000game}, which is computationally difficult for generic nonlinear systems.

The second type of methods proposed recently involves learning-based approaches. Unlike the optimal control formulation which considers the entire state space, learning-based methods rely on a finite data-set. Supervised learning-based methods have been proposed \cite{srinivasan2020synthesis,robey2020learning}, where a demonstrator is required to collect data. A neural network with a loss function encoding the conditions that a CBF needs to satisfy is used in \cite{qin2021learning,dawson2022safe,dawson2023safe,abate2021fossil}. Learning-based methods show high flexibility for nonlinear and high order systems, and are amenable to applications to high degree-of-freedom robotics. However, 
rigorous guarantees for safety and network robustness is inherently hard for these black-box methods. At the same time, the data required by the CBF network and the controller network in the training process can be difficult to obtain, as pointed out in \cite{yang2023synthesizing}.

The third type of methods involves optimization-based approaches, especially using sum-of-squares programming \cite{book3}. Barrier functions are designed for systems with input disturbances using SOS programming in \cite{prajna2007framework,prajna2006barrier}. When controllers are taken into consideration, alternating between synthesizing a controller and CBFs to solve sequential SOS programs is proposed in \cite{wang2022safety,clark2021verification,wang2018permissive,dai2022convex,kang2023verification}. Convex quadratic CBFs, constructed from a Lyapunov function for a polytopic safe set are considered in \cite{thomas2018safety,he2022barrier,zhao2023convex}. Newton's method can be leveraged to guarantee local convergence to a feasible CBF \cite{clark2022semi}. As a dual to SOS programming, moment problems based on occupation measures have been proposed \cite{korda2014convex,oustry2019inner}. These SOS-based methods transform the algebraic conditions for CBF to polynomial positivity conditions, and cast these conditions using SOS hierarchies. Compared with numerical methods to solve the Hamilton-Jacobi partial differential equations, SOS programming based methods are computationally more efficient provided that the polynomial basis is fixed. Compared with learning-based methods, SOS-based approaches allow for rigorous safety guarantee providing a feasible solution exactly. Our proposed method belongs to the SOS-based methods, whilst providing computational efficiency improvements and feasibility guarantees.

\subsection{Contribution}
In this paper, we focus on linear systems. Our main contribution is to propose an efficient method to design a control barrier function and an associated affine state feedback controller using sum-of-squares programming. The control barrier function and feedback controller are synthesized in one unified sum-of-squares program, thus overcoming the need for iterative algorithms \cite{wang2022safety,clark2021verification,wang2018permissive,dai2022convex,kang2023verification}. Moreover, our formulation is applicable to high and mixed relative degree cases without using backstepping. 

We also extend the existing literature when considering limits in the system inputs. $\mathcal{L}$-1 norm constrained limitation set is considered in \cite{wang2022safety,clark2022semi,dai2022convex,he2022barrier}. Specifically, \cite{wang2022safety,clark2022semi} introduce bilinear constraints in the sum-of-squares programming, \cite{dai2022convex} proposes a quantifier exchange to drop the dependency on the control input, and \cite{he2022barrier} proposes re-parameterization for linear systems. In our work, $\mathcal{L}$-1, $\mathcal{L}$-2, and $\mathcal{L}-\infty$ norm constrained limitations are all addressed by means of convex constraints. These input constraints can be appended to the CBF and controller synthesis program. 

\subsection{Organization}
Section \ref{sec:preliminary} provides some background. The convex synthesis program and extensions for linear systems are presented in Section \ref{sec:ConLinear}. Simulation results are shown in Section \ref{sec:simulation}. Section \ref{sec:conclusion} concludes the paper.

\section{Preliminaries}
\label{sec:preliminary}
\subsection{Notation}
For a function $b(x):\mathbb{R}^n\to\mathbb{R}$, a set denoted by the corresponding calligraphic letter $\mathcal{B}$ is defined by $\mathcal{B}:\{x\in\mathbb{R}^n:b(x)\ge 0\}$. For the set $\mathcal{B}$, $\mathcal{B}^c$ denotes the closure of its complement, $\partial \mathcal{B}$ denotes its boundary. For a positive integer $n$, $I_n$ denotes the $n\times n$ identity matrix. A positive semi-definite matrix $A$ is denoted by $A\succeq 0$. $A_{ij}$ is the element of $i$-th row and $j$-th column. For a vector $a$, $a_i$ denotes the $i$-element. $\mathrm{Tr}(A)$ is the trace of matrix $A$. $\Sigma[x]$ denotes the set of sum-of-squares polynomials in $x$. A $n \times n$ diagonal matrix is defined by $\mathrm{diag}(l_1,\ldots,l_n)$, where $l_1,\ldots,l_n\in\mathbb{R}$. For $c\in\mathbb{R}^n$ and $\epsilon>0$, $\mathcal{E}(c,\epsilon):=\{x\in\mathbb{R}^n:||x-c||_2^2\le \epsilon\}$.
\subsection{Safety and CBF}
{  
Consider a continuous-time nonlinear control-affine system
\begin{equation}\label{eq:nonlnsys}
    \dot x=f(x)+g(x)u,
\end{equation}
with $x(t) \in\mathbb{R}^n$, $u(t)\in\mathcal{U}\subseteq \mathbb{R}^m$, $f(x):\mathbb{R}^n\to\mathbb{R}^n,$ and $g(x):\mathbb{R}^n\to \mathbb{R}^{n\times m}$. Both functions are further assumed to be locally Lipschitz continuous. 
Our goal is to design a state feedback controller $u(x)$ such that the solution $x(t,x_0)$ of the closed-loop system $\dot x=f(x)+g(x)u(x)$ that starts from $x(0)=x_0$, with $x_0$ belonging to a set of \emph{initial conditions} $\mathcal{I}$, stays within a \emph{safe set} $\mathcal{S}$ for every $t$ that belongs to the domain of definition of the solution. If such a controller $u(\cdot)$ exists, we say the system is \emph{safe}. 

By a slight abuse of terminology, for our purposes we introduce the following definition:
\begin{definition} [Forward Invariance]\label{def:invariant}
    Consider system $\dot x = F(x)$, where $F:\mathbb{R}^n\to \mathbb{R}^n$ is a locally Lipschitz continuous vector field.      A set $\mathcal{C}=\{x\in \mathbb{R}^n\colon c(x)\ge 0\}$, where $c\colon \mathbb{R}^n\to \mathbb{R}$ is a continuously differentiable function,  is forward invariant for $\dot x = F(x)$ if  
\begin{equation}\label{forw.invariance}
\frac{\partial c(x)}{\partial x}F(x)\ge 0, \quad  \forall x\in \partial \mathcal{C}.
\end{equation}
If $F(x)=f(x)+g(x)u(x)$, we will refer to $\mathcal{C}$ as a {\em control} invariant set, to the function $c(\cdot)$ as a {\em control barrier function} (CBF) and to $u(x)$ as a {\em safe controller}. 
\end{definition}


Local Lipschitz continuity of $F$ implies local existence and uniqueness of the solution.  The requirement \eqref{forw.invariance} guarantees that the solution remains within  the set $\mathcal{C}$ throughout its interval of definition. 
}

\subsection{Sum-of-Squares Programming}
\begin{definition}\label{def:sos}
A polynomial $p(x)$ is said to be a sum-of-squares polynomial in $x\in\mathbb{R}^n$ if there exist {$M$} polynomials $p_i(x)$, {$i=1,\ldots,M,$} such that
\begin{equation}\label{eq:sosdecomposition}
    p(x)=\sum_{i}^Mp_i(x)^2.
\end{equation}
\end{definition}
We also call \eqref{eq:sosdecomposition} a sum-of-squares decomposition for $p(x)$. Clearly, if a function $p(x)$ has a sum-of-squares decomposition, then it is non-negative for all $x\in\mathbb{R}^n$. Computing the sum-of-squares decomposition \eqref{eq:sosdecomposition} can be efficient as it is equivalent to a positive semidefinite feasibility program.
\begin{lemma}\label{lem:sos}
Consider a polynomial $p(x)$ of degree $2d$ in $x\in\mathbb{R}^n$. Let $z(x)$ be a vector of all monomials of degree less than or equal to $d$. Then $p(x)$ admits a sum-of-squares decomposition if and only if
\begin{equation}\label{eq:sossdp}
    p(x) = z(x)^\top Q z(x), Q \succeq 0.
\end{equation}
\end{lemma}
In Lemma \ref{lem:sos}, $z(x)$ is a user-defined monomial basis if $d$ and $n$ are fixed. In the worst case, $z(x)$ has $\left( \begin{array}{c}
n + d\\
d
\end{array} \right)$ components, and $Q$ is a $\left( \begin{array}{c}
n + d\\
d
\end{array} \right) \times \left( \begin{array}{c}
n + d\\
d
\end{array} \right)$ square matrix. The necessity of Lemma \ref{lem:sos} is natural from the definition of positive semi-definite matrix, considering the monomial $z(x)$ as a vector of new variables $z_i$. The sufficiency is shown by factorizing $Q=L^\top L$. Then $z(x)^\top Q z(x)=(Lz(x))^\top Lz(x)=||Lz(x)||^2_2\ge 0$.

Given $z(x)$, finding $Q$ to decompose $f(x)$ as in \eqref{eq:sossdp} is a semi-definite program, which can be solved efficiently using interior point methods. Selecting the basis $z(x)$ depends on the structure of $p(x)$ to be decomposed. 
\begin{definition}\label{def:semi-algebrai}
A set $\mathcal{X}\subset \mathbb{R}^n$ is \emph{semi-algebraic} if it can be represented using polynomial equality and inequality constraints. If there are only equality constraints, the set is \emph{algebraic}.
\end{definition}

\begin{lemma}[S-procedure]\label{lem:S-procedure}
Suppose $t(x)\in\Sigma[x]$, then
\begin{equation}\label{eq:S-procedure1}
    p(x)-t(x)q(x)\in\Sigma[x]\Rightarrow p(x)\ge 0, \forall x\in\{x|q(x)\ge 0\}.
\end{equation}

Suppose $l(x)\in\mathbb{R}[x]$, then
\begin{equation}\label{eq:S-procedure2}
    p(x)-l(x)q(x)\in\Sigma[x]\Rightarrow p(x)\ge 0, \forall x\in\{x|q(x)=0\}.
\end{equation}
\end{lemma}
In general, compared with the \textit{Positivstellensatz}\cite{parrilo2000structured}, the S-procedure only gives a sufficient condition for the emptiness of a semi-algebraic set.
\section{Convex Design for Linear Systems}\label{sec:ConLinear}
{In this section, we propose convex synthesis programs to construct a CBF and an affine safe feedback controller. In Section \ref{sec:global-design}, we first consider a \emph{global} design for $\mathcal{B}=\{{x\in\mathbb{R}^n}:b(x)\ge 0\}$ to be control invariant. For this case, we consider the unsafe set $\mathcal{S}^c$ to be bounded on a subspace of $\mathbb{R}^n$. This is commonly for robot collision avoidance problems, where the position space is a subspace of the robot state space. The control invariant set $\mathcal{B}$ is constructed \emph{globally} as its projection to the subspace of $\mathcal{S}^c$ is unbounded. For the second case in Section \ref{sec:local-design}, we construct a control invariant set $\mathcal{B}^c=\{{x\in\mathbb{R}^n}:b(x)\le 0\}$ around a bounded initial set. This control invariant set is called \emph{local} as we will show it is bounded on $\mathbb{R}^n$.} 

{ 
Consider a continuous-time linear system:
\begin{equation}\label{eq:linearsys}
    \dot x = Ax+Bu,
\end{equation}
where $x(t)\in\mathbb{R}^n$, $u(t)\in\mathcal{U}\subseteq \mathbb{R}^m$ are the state and control input, and $A\in\mathbb{R}^{n\times n}$, $B\in\mathbb{R}^{n\times m}$. We assume that the system is stabilizable. {Throughout the paper, the CBF $b(x)$ and feedback controller $u(x)$ are parameterized as follows.
\begin{subequations}
    \begin{align}
    b(x)&=(x-c)^\top \Omega^{-1}(x-c)-1,\label{eq:cbf}\\
    u(x)&=Y\Omega^{-1}(x-c)+d\label{eq:control},
\end{align}
\end{subequations}
where $\Omega\in\mathbb{R}^{n\times n},Y\in\mathbb{R}^{m\times n}$ are matrices to be designed, $c\in\mathbb{R}^n$ and $d\in\mathbb{R}^m$ are constant vectors that will be clear in the sequel.}
}


\subsection{Global Design}\label{sec:global-design}
Consider the safe set $\mathcal{S}$ defined by a union of semi-algebraic sets as:
\begin{equation}\label{eq:globalsafe}
    \mathcal{S}:=\bigcup_{i=1}^o\{x\in\mathbb{R}^n:s_i(\overline x)\ge 0\},
\end{equation}
where $x=[\overline x^\top,\underline x^\top]^\top$, with 
$\overline x\in\mathbb{R}^{\overline n}$, $\underline x\in\mathbb{R}^{\underline n}$ and $\overline n+\underline n=n$.
If $\underline n=0$, the safe set is defined over all the states.

\begin{assumption}\label{ass:safeglobal}
    $\mathcal{S}$ is a semi-algebraic set, and $\mathcal{S}^c$ is bounded on the space $\mathbb{R}^{\overline n}$.
\end{assumption}
Let $c=[\overline c^\top,\underline c^\top]^\top\in\mathbb{R}^{n}$ be a vector of constants such that rank($[B,Ac]$) = rank($B$), and consider the following optimization {program}:
\begin{subequations}\label{eq:dqcbflem}
\begin{align}
\min~&\mathrm{Tr}(\overline{\Omega})  \\
\mathrm{subject~to}~ &{0}\prec \overline\Omega = \overline\Omega^\top \in\mathbb{R}^{\overline n\times \overline n}, {0}\succ\underline \Omega = \underline \Omega^\top \in\mathbb{R}^{\underline n\times \underline n}{}, \label{eq:dqcbflemeq1}\\
&{0}\prec R = R^\top \in\mathbb{R}^{\overline n\times \overline n}\succ 0, Y\in\mathbb{R}^{m\times n}, \\
&\sigma_1(\overline x),\ldots,\sigma_o(\overline x)\in\Sigma[\overline x], \epsilon>0\\
&\Omega = \left[ {\begin{array}{*{20}{c}}
{{\overline \Omega}}&{0}\\
{0}&{{\underline \Omega}}
\end{array}} \right] \label{eq:cascade} \\
& \Omega A^\top+Y^\top B^\top +A\Omega + BY \succeq 0 \label{eq:dqcbflemeqa} \\
&\begin{bmatrix}
            R&I_{\overline n}\\I_{\overline n}&\overline\Omega 
        \end{bmatrix}\succeq 0 \label{eq:dqcbflemeqb} \\
 & 1-\overline x_c^\top R\overline x_c+\sum_{i=1}^o\sigma_i(\overline x)s_i(\overline x)-\epsilon\in\Sigma[\overline x], \label{eq:dqcbflemeqc}  
\end{align}  
\end{subequations}
where $\overline x_c=\overline x-\overline c$, $\underline x_c=\underline x-\underline c$, $x_c=x-c$. Notice 
that this a convex optimization program, where the objective function is linear, and is subject to semi-define constraints
{The cost function is to minimize the volume of the set $\{\overline{x}\in\mathbb{R}^{\overline{n}}:-(\overline{x}-\bar c)^\top \overline{\Omega}^{-1}(\overline{x}-\bar c)+1\ge 0\}$, thus indirectly maximizing the volume of the projection set of $\mathcal{B}$ on the space $\mathbb{R}^{\overline n}$. An alternative formulation is $\max\log\det\overline{\Omega}^{-1}$ \cite[Section 2.2.4]{boyd1994linear}. However, this is not supported by SeDuMi, which is the solver we are using to solve the semi-definite program.}
In the following theorem, we give the main result of the paper, a convex program to synthesize a CBF $b(x)$ and a feedback controller $u(x)$ under Assumption \ref{ass:safeglobal}. 

\begin{theorem}\label{th:dqcbf}
Consider Assumption \ref{ass:safeglobal}, and let $\mathcal{U}=\mathbb{R}^m$. Assume that a solution to \eqref{eq:dqcbflem} exists and is denoted by $\overline\Omega, \underline \Omega, R, Y, \{\sigma_i(\cdot)\}_{i=1}^o, \epsilon$. 
Set $u(x)=Y\Omega^{-1}{(x-c)}+d$ where
$d\in\mathbb{R}^m$ is such that $Bd+Ac = 0$. 
We then have that 
{ 
\begin{enumerate}
\item $\mathcal{B}:=\{x\in\mathbb{R}^n:{(x-c)}^\top \Omega^{-1}{(x-c)}-1\ge 0\}\subseteq\mathcal{S}$.
\item $\mathcal{B}$ is a control invariant set for $\dot x = Ax+Bu(x)$.
\end{enumerate}
}
\end{theorem}

\begin{proof}
We first prove that satisfaction of 
\eqref{eq:cascade}, \eqref{eq:dqcbflemeqb} and \eqref{eq:dqcbflemeqc} are sufficient for $\mathcal{B}\subseteq \mathcal{S}$.
Given that $R\succ 0$ and $\overline\Omega \succ 0$, using Schur complement, \eqref{eq:dqcbflemeqb} is equivalent to $R-\overline\Omega^{-1}\succeq 0$. Multiplying the latter condition by {$\overline{x}_c^\top=(\overline{x}-\overline{c})^\top$} on the left and by $\overline{x}_c$ on the right, we obtain 
\begin{equation*}
    -\left(\overline x_c^\top \overline\Omega^{-1}\overline x_c-1\right)+(\overline x_c^\top R \overline x_c-1)\ge 0, \forall \overline x\in\mathbb{R}^{\overline n}.
\end{equation*}
Then we have the following relationship
\begin{align*}
    \{\overline{x}\in\mathbb{R}^{\overline{n}}&:{(\overline{x}-\overline{c})^\top} \overline{\Omega}^{-1}{(\overline{x}-\overline{c})}-1> 0\} \nonumber \\
    &\subseteq \{\overline x\in\mathbb{R}^{\overline{n}}:{(\overline{x}-\overline{c})}^\top R{(\overline{x}-\overline{c})}-1> 0\}. 
\end{align*}
{The former set inclusion implies the following relationship for the closures of the associated sets, 
\begin{align*}
    \{\overline{x}\in\mathbb{R}^{\overline{n}}&:(\overline{x}-\overline{c})^\top \overline{\Omega}^{-1}(\overline{x}-\overline{c})-1\ge 0\} \nonumber \\
    &\subseteq \{\overline x\in\mathbb{R}^{\overline{n}}:(\overline{x}-\overline{c})^\top R(\overline{x}-\overline{c})-1\ge 0\}.
\end{align*}
The two sets involved are subsets of $\mathbb{R}^{\overline n}$. Considering them as the base of cylinder sets in $\mathbb{R}^n$, we obtain
\begin{align}
    \{{x}\in\mathbb{R}^{{n}}&:(\overline{x}-\overline{c})^\top \overline{\Omega}^{-1}(\overline{x}-\overline{c})-1\ge 0\} \nonumber \\
    &\subseteq \{ x\in\mathbb{R}^{{n}}:(\overline{x}-\overline{c})^\top R(\overline{x}-\overline{c})-1\ge 0\}.
    \label{eq:bsuff1}
\end{align}
}
Invoking Lemma \ref{lem:S-procedure} for polynomial functions $s_i(\overline x)$, $i=1,\ldots,o$, \eqref{eq:dqcbflemeqc} indicates that
\begin{equation}
    -{(\overline{x}-\overline{c})}^\top R{(\overline{x}-\overline{c})}+1>0~\forall x:s_1(\overline x),\ldots,s_o(\overline x)<0,
\end{equation}
{
which in turn implies
\begin{align}
    \exists i\in\{1,2,\ldots,o\}&:s_i(\overline{x})\ge0,~\forall {x}\in\mathbb{R}^{ n}~\nonumber \\
    &\mathrm{such~that}~(\overline{x}-\overline{c})^\top R(\overline{x}-\overline{c})-1\ge 0.
\end{align}}
Converting the above relationship in set inclusion form, we have
{
\begin{equation}\label{eq:bsuff2}
    \{{x}\in\mathbb{R}^{{n}}:(\overline{x}-\overline{c})^\top R(\overline{x}-\overline{c})-1\ge 0\}\subseteq \mathcal{S}.
\end{equation}}
Combining \eqref{eq:bsuff1} with \eqref{eq:bsuff2} implies 
\begin{equation}\label{eq:21}
    \{{x}\in\mathbb{R}^{{n}}:{(\overline{x}-\overline{c})}^\top \overline{\Omega}^{-1}{(\overline{x}-\overline{c})}-1\ge 0\}\subseteq \mathcal{S}.
\end{equation}
Given that $\underline \Omega \prec 0$, and using the block representation in \eqref{eq:cascade}, we have 
\begin{align}
    {(x-c)}^\top\Omega^{-1} {(x-c)}-1&={(\overline{x}-\overline{c})}^\top\overline\Omega^{-1}{(\overline{x}-\overline{c})}\nonumber\\+{(\underline{x}-\underline{c})}^\top\underline\Omega^{-1}{(\underline{x}-\underline{c})}&-1 \le {(\overline{x}-\overline{c})}^\top\overline\Omega^{-1}{(\overline{x}-\overline{c})}-1. \label{eq:22}
\end{align}
Using \eqref{eq:22} into the relationship in \eqref{eq:21}, and recalling that $\mathcal{B}=\{x\in\mathbb{R}^n:{(x-c)^\top}\Omega^{-1}{(x-c)} -1 \ge 0\}$, we obtain
\begin{equation}
     \mathcal{B}\subseteq \{{x}\in\mathbb{R}^{{n}}:{(\overline{x}-\overline{c})}^\top \overline{\Omega}^{-1}{(\overline{x}-\overline{c})}-1\ge 0\}\subseteq \mathcal{S}.
\end{equation}

We then prove that $\mathcal{B}$ (the zero super-level set of $b(x)$) is a control invariant set. Since $\Omega$ (decomposed as in \eqref{eq:cascade}), is optimal for \eqref{eq:dqcbflem}, it will have to satisfy \eqref{eq:dqcbflemeqa}. We will show that \eqref{eq:dqcbflemeqa} is sufficient for $\mathcal{B}$ to be control invariant, thus establishing the claim.
We guarantee control invariance by using an affine state feedback controller as \eqref{eq:control},
where $K\in\mathbb{R}^{m\times n}$.
Transforming the coordinate from $x$ to {$x_c=x-c$}, and since $d$ is such that $Bd +Ac = 0$, the transformed system dynamics are given by
\begin{equation*}
    \dot{x}_c = (A+BK)x_c.
\end{equation*}
In the new coordinate $b(x)|_{x=x_c+c}=\tilde b(x_c)=x_c^\top \Omega^{-1}x_c-1$.
If $\dot b(x)\ge 0$ for any $x\in\mathbb{R}^n$, then $\mathcal{B}$ is invariant. Notice that $\dot b(x)\ge 0$ is equivalent to
\begin{align}
    \dot b(x)=\dot{\tilde b}(x_c) = &~x_c^\top(A^\top\Omega^{-1}+\Omega^{-1} A\nonumber \\
    &+K^\top B^\top\Omega^{-1}+\Omega^{-1}BK)x_c\ge 0. \label{eq:evaluationofbx}
\end{align}
Satisfaction of \eqref{eq:evaluationofbx} for any $x\in\mathbb{R}^n$ is equivalent to
\begin{align}\label{eq:intereq}
    &A^\top\Omega^{-1}+\Omega^{-1} A+K^\top B^\top\Omega^{-1}+\Omega^{-1}BK\succeq 0.
\end{align}
Left and right multiplying by $\Omega$ on both sides of \eqref{eq:intereq}, we obtain
\begin{align*}
    &\Omega A^\top + A\Omega+\Omega K^\top B^\top +BK\Omega\succeq 0.
\end{align*}
Substituting $K\Omega$ with a new matrix $Y\in\mathbb{R}^{m\times n}$, and noticing that $\Omega$ is invertible, we equivalently obtain \eqref{eq:dqcbflemeqa}. The latter is thus a sufficient condition for 
$\mathcal{B}$ to be invariant.

By the proof of the previous part it follows that $u(x)=K{(x-c)}+d$ renders $\mathcal{B}$ invariant. Moreover, $Y = K \Omega$, which in turn implies that $K = Y \Omega^{-1}$. Therefore, $u(x)=Y\Omega^{-1}{(x-c)}+d$ guarantees invariance. 
\end{proof}
\begin{figure}[t]
    \centering
    \includegraphics[scale=0.3]{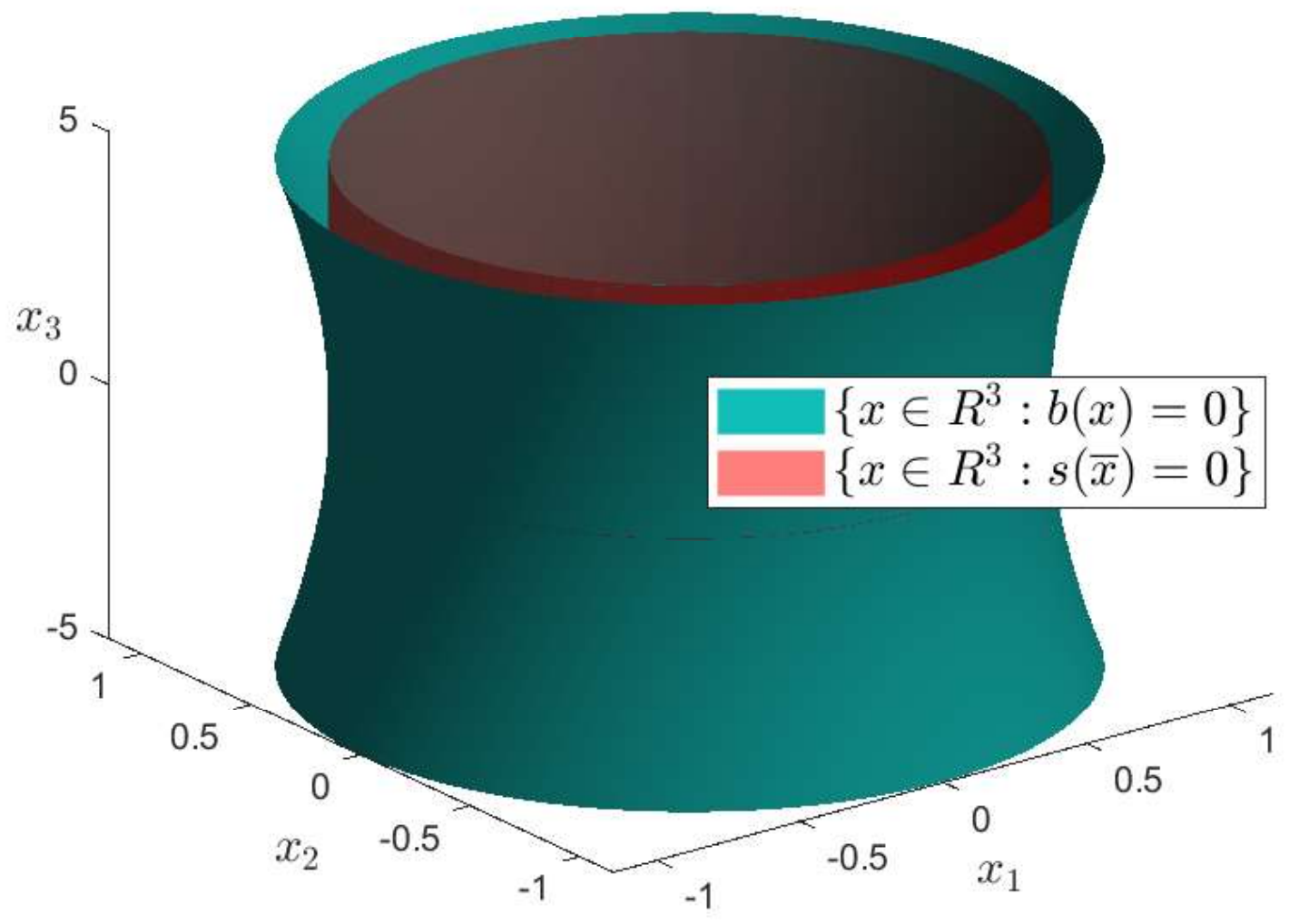}
    \caption{{ {Visualization of Case \ref{cas:2} in Section \ref{sec:case}}, where the state {$\overline x=[{x}_1,{x}_2]\in\mathbb{R}^2$, $\underline{x}=x_3\in\mathbb{R}$}. The safe set $\mathcal{S}:=\{x\in\mathbb{R}^3:s(\overline{x})\ge 0\}$ is a cylinder expanded from a set on $\mathbb{R}^2$ to $\mathbb{R}^3$. The region outside of the blue hyperboloid represents the set $\mathcal{B}:=\{x\in\mathbb{R}^3:b(x)\ge 0\}$, which is control invariant from our construction. The safe set $\mathcal{S}:=\{x\in\mathbb{R}^3:s(\overline{x})\ge 0\}$ is the outside of the {inner} red cylinder. We can see from the figure that $\mathcal{B}\subseteq \mathcal{S}$.}}
    \label{fig:th1}
\end{figure}
Figure \ref{fig:th1} visualizes an example of a control invariant set $\mathcal{B}$ and $\mathcal{S}$ on $\mathbb{R}^3$.
{To gain an intuitive understanding of Theorem \ref{th:dqcbf}, we analytically construct a control invariant set for a simple example.}

{ 
\begin{example}\label{ex:planar}
{Consider a car moving along a line}
\begin{align*}
    &\dot{\overline{x}} = \underline{x}\\
    &\dot{\underline{x}} = u
\end{align*}
where $\overline{x}$ represents the position and $\underline{x}\in\mathbb{R}$ represents the velocity. Let
\[
\mathcal{S} = \{x\in\mathbb{R}^2:s(\overline{x}):=\overline{x}^2-1\ge 0\}.
\]
We follow the construction in Theorem \ref{th:dqcbf}.
The vector $c\in\mathbb{R}^2$ that satisfies rank($[B~Ac]$) = rank($B$) is any vector such that $\underline{c}=0$. We also fix $\overline{c} = 0$. Consequently $d=0$. 

Consider the decision variables $\overline{\Omega}\in{\mathbb{R}_{>0}}$, $\underline{\Omega}\in{\mathbb{R}_{<0}}$, $R\in\mathbb{R}$, $Y=[Y_1~Y_2]\in\mathbb{R}^{1\times2}$, $\sigma(\overline{x})\in\Sigma[\overline{x}]$ and $\varepsilon>0$. The constraints in \eqref{eq:dqcbflem} can be written as follows
\begin{subequations}
    \begin{align}
        \Omega = \begin{bmatrix}
                \overline \Omega & 0 \\
                0 & \underline \Omega
            \end{bmatrix},\label{eq:ex1eq1}\\
        \begin{bmatrix}
            0 & \underline\Omega+Y_1\\
            \underline\Omega+Y_1&2Y_2
        \end{bmatrix}\succeq 0,\label{eq:ex1eq2}\\
        \begin{bmatrix}
            R&1\\1&\overline{\Omega}
        \end{bmatrix}\succeq 0,\label{eq:ex1eq3}\\
        1-R\overline{x}^2+\sigma(\overline{x})(\overline{x}^2-1)-\epsilon\in\Sigma[\overline{x}].  \label{eq:ex1eq4}\end{align}
\end{subequations}
Looking at the spectrum of the matrix in \eqref{eq:ex1eq2}, we conclude that \eqref{eq:ex1eq2} is equivalent to $Y_2>0$ and $\underline{\Omega}+Y_1=0$. Condition \eqref{eq:ex1eq3} is equivalently expressed as $R-\overline{\Omega}^{-1}\ge 0$. As for $\sigma(\overline{x})$, we set it to be an SOS polynomial of degree $0$, hence, $\sigma(\overline{x})=\overline{\sigma}\ge 0$. Writing the polynomial in \eqref{eq:ex1eq4} as
\begin{equation*}
    1-R\overline{x}^2+\sigma(\overline{x})(\overline{x}^2-1)-\epsilon=\begin{bmatrix}
        1\\\overline{x}^2
    \end{bmatrix}^\top \begin{bmatrix}
        1-\overline{\sigma}-\epsilon&0\\
        0&-R+\overline{\sigma}
    \end{bmatrix}
    \begin{bmatrix}
        1\\\overline{x}^2
    \end{bmatrix}
\end{equation*}
and bearing in mind Lemma \ref{lem:sos}, one realizes that condition \eqref{eq:ex1eq4} is equivalent to ${\bar\sigma}-R\ge 0$ and ${1-}\overline{\sigma}-\epsilon\ge 0$. In summary, we have the following conditions
\begin{align*}
    &Y_2>0~\mathrm{and}~Y_1=-\underline{\Omega}>0\\
    &\overline{\sigma}\ge R\ge \overline{\Omega}^{-1}>0\\
    &\epsilon>0~\mathrm{and}~\overline{\sigma}+\epsilon\le 1
\end{align*}
We set $Y=[Y_1~Y_2]:=[-\underline \Omega~Y_2]$, with $\underline{\Omega},Y_2$ any positive numbers, and $\overline{\Omega}>1$, $R=\overline{\sigma}=\overline{\Omega}^{-1}$, $\epsilon\le 1-\overline{\sigma}$. Note that setting $\overline{\Omega}>1$ is necessary for having the constraint $\overline{\sigma}+\epsilon\le 1$ satisfied. To minimize the cost function, we should take $\overline{\Omega}$ as {small} as possible. We obtain that the function $b(x)$ that defines $\mathcal{B}:=\{x\in\mathbb{R}^n:b(x)\ge 0\}$ is
\[
b(x)=\overline{\Omega}^{-1}\overline{x}^2-\underline\Omega^{-1}\underline{x}^2-1,
\]
and the feedback controller is
\[
u = Y\Omega^{-1}=[-\underline{\Omega}\slash\overline\Omega~Y_2\slash\underline\Omega],
\]
resulting in the closed-loop dynamics
\[
\dot x = \begin{bmatrix}
    0&1\\
    -\underline\Omega\slash\overline{\Omega}&Y_2\slash\underline\Omega
\end{bmatrix}x
\]
Figure \ref{fig:ex1} shows a sketch of the {sets} $\mathcal{S}^c$, $\mathcal{B}$. We first observe that $\overline{\Omega}>1$ established above guarantees that $\mathcal{B}\subset \mathcal{S}$. Second, formulating $\mathcal{B}$ in terms of the entire state vector $x$ results in a set $\mathcal{B}$ which differs from the one a designer could expect, namely, $\mathcal{B}:=\{x\in\mathbb{R}^n:{|\overline{x}|}\ge \overline{b}\}$, with $\overline{b}>1$.

For any feasible choice of the design parameters, the obtained closed-loop matrix has at least one unstable eigenvalue. To have an understanding of the state response, we compute the spectral representation $e^{(A+BK)t}$ for these values of the design parameters: $\underline\Omega=-4$, $\overline{\Omega}=2$, $Y_2=4$. 
Then the spectral representation is given by
\begin{align*}
    e^{(A+BK)t}=\frac{1}{3}\begin{bmatrix}
        e^{-2t}+2e^t&-e^{-2t}+e^t\\
        -2e^{-2t}+2e^t&2e^{-2t}+e^t
    \end{bmatrix}.
\end{align*}
Hence, if the system starts from the initial condition $x=\begin{bmatrix}
    \overline{\Omega}&0
\end{bmatrix}^\top = [2~0]^\top $, which is on the boundary of $\mathcal{B}$, it will evolve as
\[
x(t)=\frac{2}{3}\begin{bmatrix}
    e^{-2t}+2e^t\\-2e^{-2t}+2e^t
\end{bmatrix}.
\]
As a result, both position and velocity diverge exponentially but are certified to stay within $\mathcal{B}$.
    
\end{example}

The CBF $b(x)$ constructed in the {previous example} is a function of the whole state $x=[\overline{x},\underline{x}]^\top$. {However, given the definition of $\mathcal{S}^c$, which only constrains the position variable $\overline x$, one could alternatively consider a candidate barrier function $\overline{b}(x):=(\overline{x}-\overline{c})^\top \overline{\Omega}^{-1}(\overline{x}-\overline{c})-1$ and the corresponding set $\overline{\mathcal{B}}:=\{x\in\mathbb{R}^n:\overline{b}(x)\ge 0\}$.} We will show {below} that the set {$\overline{\mathcal{B}}$} can not be control invariant using linear feedback $u(x)$. Denote the projection matrix
\[\Pi:=[I_{\overline{n}}~~0_{\overline{n}\times {\underline n}}].\]
{Use $\overline{x}_c$ instead of $\overline{x}-\overline{c}$, and let $\overline{x}_c=\Pi x_c$.} We can derive the following {identity}
\[
\dot {\overline{b}}(x)=x_c^\top(\Pi^\top \overline{\Omega}^{-1}\Pi(A+BK)+(A+BK)^\top \Pi^\top \overline{\Omega}^{-1}\Pi)x_c
\]
and express the invariance condition as
\[
\Pi^\top \overline{\Omega}^{-1}\Pi(A+BK)+(A+BK)^\top \Pi^\top \overline{\Omega}^{-1}\Pi\succeq 0.
\]
\begin{figure}[t]
    \centering
    \includegraphics[scale=0.3]{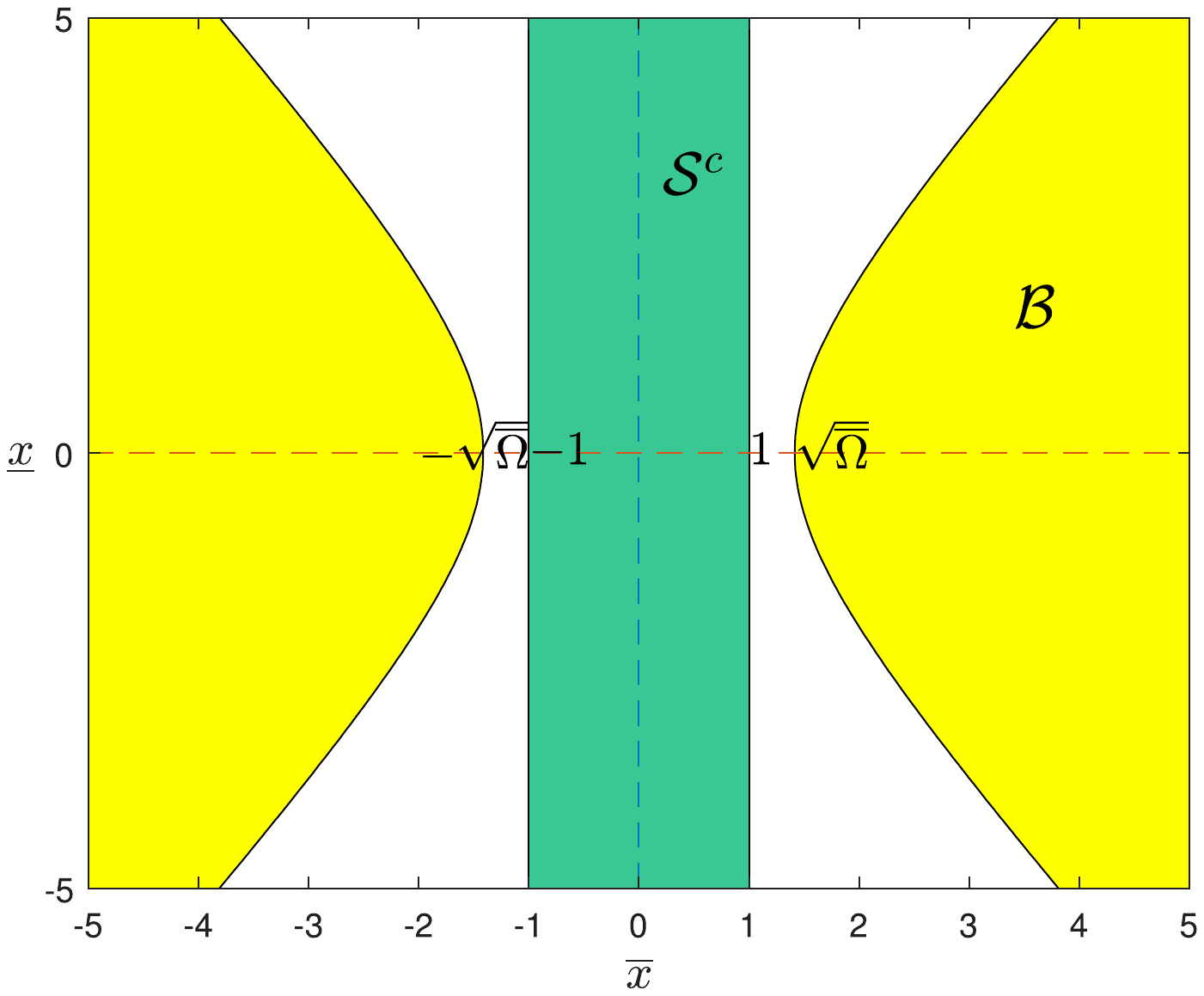}
    \caption{{  Pictorial illustration for Example \ref{ex:planar}. The green set $\mathcal{S}^c:=\{x\in\mathbb{R}^2:\overline{x}^2-1{\le} 0\}$ is expanded from a segment on $\mathbb{R}^1$. The designed control invariant set $\mathcal{B}:=\{x\in\mathbb{R}^2:\overline{\Omega}^{-1}\overline{x}^2-\underline{\Omega}\underline{x}^2-1\ge 0$\} has been filled in yellow. Intuitively, with a large velocity, i.e. larger $|\underline x|$, the planar car should stay further away from the obstacle, which can be seen by the gap between $\mathcal{B}$ and $\mathcal{S}^c$ {being larger for larger $|\overline{x}|$.}}}
    \label{fig:ex1}
\end{figure}
We partition $A+BK$ according to the partition $\mathbb{R}^{\overline{n}}\times \mathbb{R}^{\underline{n}}$ to obtain
\[
A+BK=\begin{bmatrix}
    \overline{A}_1&\overline{A}_2\\
    \underline{A}_1&\underline{A}_2
\end{bmatrix}+
\begin{bmatrix}
    \overline{B}\\\underline{B}
\end{bmatrix}
\begin{bmatrix}
    \overline{K}&\underline{K}
\end{bmatrix}
\]
The invariance condition can be expressed as
\begin{align*}
    &\begin{bmatrix}
    \overline{\Omega}^{-1}(\overline{A}_1+\overline{B}\overline{K})&\overline{\Omega}^{-1}(\overline{A}_2+\overline{B}\underline{K})\\
    0&0
\end{bmatrix}+\\
&\begin{bmatrix}
    (\overline{A}_1+\overline{B}\overline{K})^\top \overline{\Omega}^{-1}&0\\
    (\overline{A}_2+\overline{B}\underline{K})^\top \overline{\Omega}^{-1}&0
\end{bmatrix}\succeq 0,
\end{align*}
which leads to a convex condition by multiplying {$\begin{bmatrix}
    \overline{\Omega}&0\\0&I_n
\end{bmatrix}$} on both sides of the matrices in the inequality. However, the possibility of fulfilling such constraint appears to be related to the possibility of shaping the spectra of $\overline{A}_1+\overline{B}\overline{K}$ and $\overline{A}_2+\overline{B}\underline{K}$, hence, to the controllability of the pairs $(\overline{A}_1,\overline{B})$, $(\overline{A}_2,\overline{B})$. Going back to Example \ref{ex:planar}, we have $(\overline{A}_1,\overline{B})=(0,0)$, and $(\overline{A}_2,\overline{B})=(1,0)$, {which shows lack of controllability of both pairs.} As a result, the invariance condition above is
\[
\begin{bmatrix}
    0&\overline{\Omega}^{-1}\\
    \overline{\Omega}^{-1}&0
\end{bmatrix}\succeq 0
\]
which shows that enforcing invariance for the set $\overline{\mathcal{B}}$ via feedback is impossible due to the lack of controllability (the matrix has a positive and a negative eigenvalue).
\begin{figure}[t]
    \centering
    \includegraphics[scale=0.5]{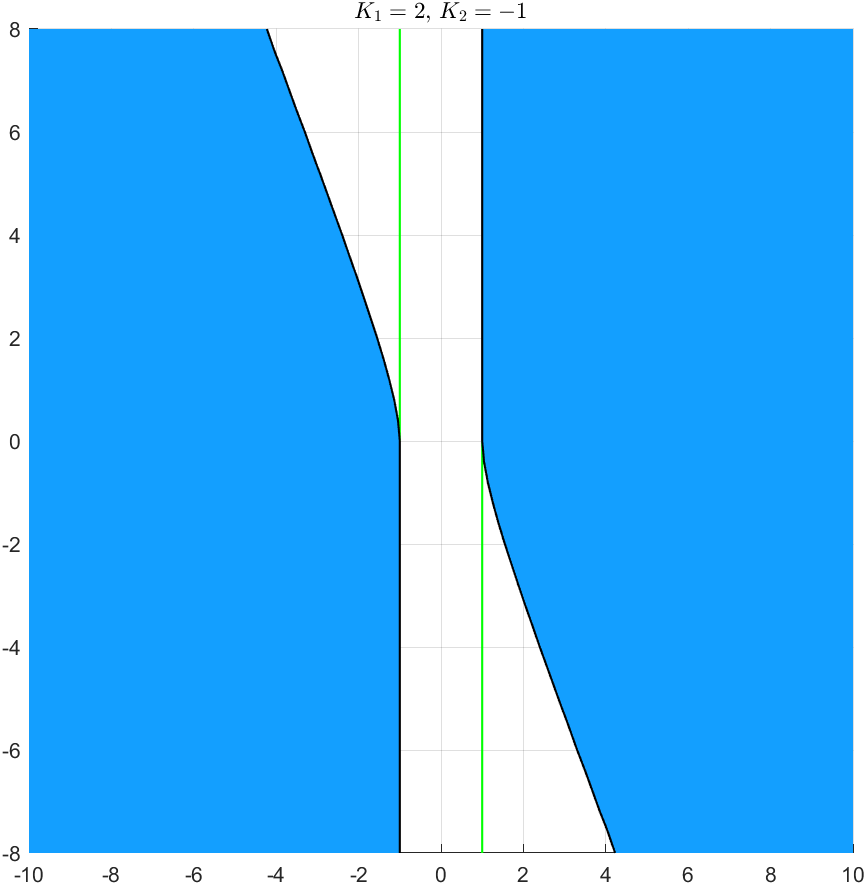}
    \caption{{ Exact invariant set for the planar car using $u=K_1\overline{x}+K_2\underline{x}$, where $K_1=Y_1\overline{\Omega}^{-1}={2},K_2=Y_2\underline{\Omega}^{-1}={-1}$. The green vertical lines are the boundary of $\mathcal{S}^c$, the set filled in blue is {the exact invariant set}.}}
    \label{fig:exact-invariant}
\end{figure}

\begin{figure}[t]
    \centering
    \includegraphics[scale=0.3]{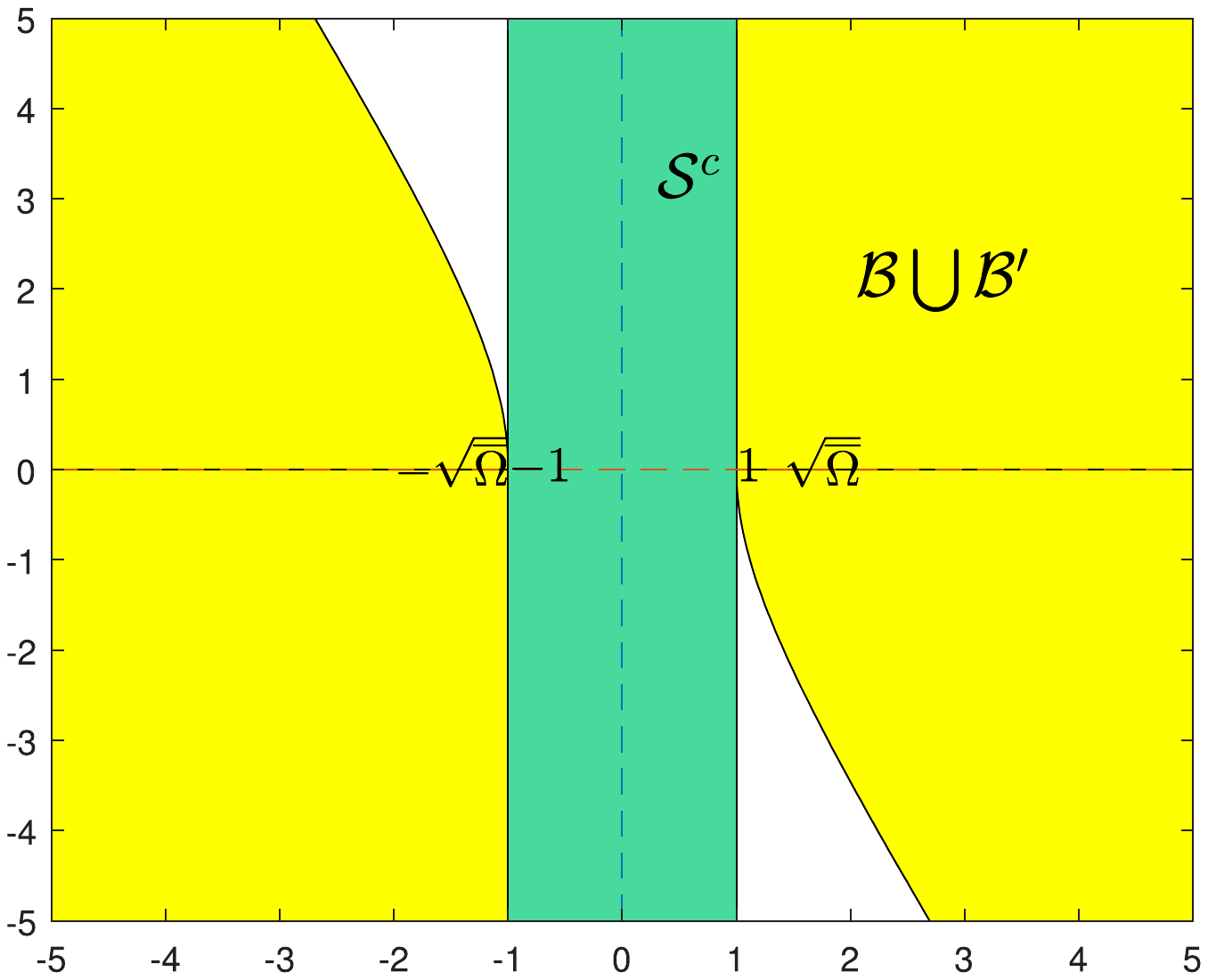}
    \caption{{The union of control invariant set $\mathcal{B}\bigcup\mathcal{B}'$. $\mathcal{B}$ is computed by solving our program \eqref{eq:dqcbflem} which results in $\overline{\Omega}=1$. Physical considerations for $\mathcal{B}'$ is obtained  from the planar car. The union is also control invariant \cite{blanchini1999set}, and close to the exact invariant set as in Figure \ref{fig:exact-invariant}.}}
    \label{fig:union-invariant}
\end{figure}
The exact control invariant set {for Example \ref{ex:planar} numerically} computed by the level-set method toolbox \cite{mitchell2007toolbox} is shown in Figure \ref{fig:exact-invariant}. {To compute the exact control invariant set, the control feedback $u$ is set in the linear} form {$u=Y\Omega^{-1}x$, where $Y$ and $\Omega$ have the same numerical values as those chosen in Example \ref{ex:planar} ($\underline\Omega=-4$, $\overline{\Omega}=2$, $Y_1=4$, $Y_2=4$).} In comparison, our computed control invariant set $\mathcal{B}$ determined analytically and depicted in Figure \ref{fig:ex1} is conservative when $\overline{\Omega}\ne 1$. This can be alleviated by minimizing $\mathrm{Tr}(\overline{\Omega})$ as in the program \eqref{eq:dqcbflem}. Conservative behaviour is also encountered in the first and third quadrants, where the boundary of the exact invariant set coincides with the safe set $\mathcal{S}$. This is natural as the planar car is moving away from the {unsafe set $\mathcal{S}^c$, which has been filled in green,} in these regions. Our method, however, computes a control invariant set $\mathcal{B}$, that is symmetric with respect to the {$\overline{x}$}-axis. In practice, one can reduce this conservativeness by taking the union of our computed control invariant set $\mathcal{B}$ with other invariant sets, such as $\mathcal{B}'=\{x\in\mathbb{R}^2:b'(x):=\underline{x}\overline{x}\ge 0\}$. The new control invariant set is shown in Figure \ref{fig:union-invariant}. The control barrier function corresponds to this union set can be defined by 
\[
\beta(x)=\max\{b(x),b'(x)\}.
\]
{Such a kind of CBF has been investigated in \cite{glotfelter2017nonsmooth}.}

According to Nagumo's Theorem \cite{nagumo1942lage}, a compact set is invariant for a vector field if and only if the vector field is within the tangent cone for all points on the boundary of the set. For a compact and closed set $\mathcal{B}$, this is equivalent to having $\dot b(x)\ge 0$, for any $x$ such that $b(x)=0$. However in our proposed convex conditions \eqref{eq:dqcbflem}, we enforce a ``strengthened" condition that $\dot b(x)\ge 0$, for any $x\in\mathbb{R}^n$. Nevertheless, we show in the following proposition that this does not introduce any conservativeness {  in the case $\overline{n}=n$}.
\begin{proposition}\label{pro:4}
    {Consider the system $\eqref{eq:linearsys}$,  constant $c\in\mathbb{R}^{n}$ such that ${\rm rank} ([B\, Ac])={\rm rank} (B)$ and a quadratic 
    function $b(x)=(x-c)^\top \Omega^{-1} (x-c)-1$, with $\Omega\succ 0$. 
    If there exists a feedback controller $u=K(x-c)+d$, with $d$ satisfying $Bd+Ac=0$, such that $\dot b(x)\ge 0$ for any $x$ such that $b(x)=0$, then $\dot b(x)\ge 0$ for any $x\in\mathbb{R}^n$.}
\end{proposition}
\begin{proof}
  {For $x=c$, $\dot b(x)= 2(x-c)^\top \Omega^{-1} (A+BK) (x-c)=0$. 
On the other hand, observe that for any point $x\ne c\in\mathbb{R}^n$, there exists $y=\frac{1}{\lambda} (x-c) +c$ with $\lambda = ((x-c)^\top\Omega^{-1}(x-c))^{1/2}>0$, such that $b(y)=0$. The function $\dot b(x)$ can be rewritten as 
\[
\dot b(x) := 2\lambda^2 \frac{(x-c)^\top}{\lambda}\Omega^{-1}(A+BK) \frac{(x-c)}{\lambda}=\lambda^2\dot b(y)
\]
As $b(y)=0$, we have $\dot b(y)\ge 0$ by the proposition’s statement that assumes this is the case for $y$ such that $b(y)=0$, which implies $\dot b(x)\ge 0$, as claimed.}
\end{proof}

As a result of Proposition \ref{pro:4}, the synthesized controller $u(x)$ endows robustness as $\dot b(x)\ge 0$ for any $x$ such that $b(x)<0$. If the system starts from an unsafe point $x$, our synthesized CBF guarantees that there exists a controller that {forces the state of the system} enters the safe region, if the problem is feasible. This property is especially helpful for unexpected perturbations to the system.


{
\begin{corollary}\label{co:1}Assume that the projection of $\mathcal{S}^c$ onto $\mathbb{R}^{\overline{n}}$ is a polytope on the space $\mathbb{R}^{\overline n}$ with vertices denoted by $v_1,\ldots,v_{o'}\in\mathbb{R}^{\overline n}$. Constraint \eqref{eq:dqcbflemeqc} can be replaced by linear constraints:
    \begin{equation}\label{eq:co1}
        -(v_i-\overline c)^\top R(v_i-\overline c)+1\ge 0,i=1,\ldots,o'.
    \end{equation}
\end{corollary}
\begin{proof}
    Constraint \eqref{eq:dqcbflemeqc} implies $\mathcal{S}^c\subseteq \mathcal{R}:=\{x\in\mathbb{R}^n:-(\overline x-\overline c)^\top R(\overline x-\overline c)+1\ge 0\}$. {Denote the projection set of }$\mathcal{S}^c$ {onto $\mathbb{R}^{\overline{n}}$ by $\overline{\mathcal{S}^c}$, which is a polytope, and} {the projection set of} $\mathcal{R}$ {onto $\mathbb{R}^{\overline{n}}$ by $\overline{\mathcal{R}}$, which is} {an ellipsoid}. We then have $\mathcal{S}^c\subseteq\mathcal{R}$ is equivalent to {$\overline{\mathcal{S}^c}\subseteq \overline{\mathcal{R}}$, which can be verified by} the constraints of all the vertices of {$\overline{\mathcal{S}^c}$} be within {$\overline{\mathcal{R}}$}. We conclude the proof.
\end{proof}
}
The number of linear constraints \eqref{eq:co1} depends on $o'$. If the polytopic $\mathcal{I}$ has $l$ facets, then the maximum number of vertices is $\left( \begin{array}{l}
o - \left\lceil {n/2} \right\rceil \\
\left\lfloor {n/2} \right\rfloor 
\end{array} \right) + \left( \begin{array}{l}
o - \left\lfloor {n/2} \right\rfloor  - 1\\
\left\lceil {n/2} \right\rceil  - 1
\end{array} \right)$, which could be quite large. For practical purposes, Corollary \ref{co:1} becomes useful if the number of vertices is moderate. 

\subsection{Local Design}\label{sec:local-design}
{In the previous section, we construct a control invariant set $\mathcal{B}:=\{x\in\mathbb{R}^n:b(x)\ge 0\}$ globally, it is unbounded on $\mathbb{R}^{ \overline{
n}}$, and naturally unbounded on $\mathbb{R}^{n}$. As shown in Example \ref{ex:planar}, the closed-loop trajectory diverges using the co-designed linear feedback controller $u(x)$. This is undesired in many applications where boundedness of trajectories is a prerequisite. In this section, we consider constructing a bounded control invariant set around a bounded {set of initial conditions} $\mathcal{I}$, and inside a intersection of half planes, i.e. the safe set $\mathcal{S}$. The new control invariant set will also be parameterized by a quadratic function. To ease notation, we still use $b(x)=x_c^\top \Omega^{-1}x_c-1$, but the new control invariant set will be derived by a sub-level set of the function, i.e. $\mathcal{B}^c:=\{x\in\mathbb{R}^n:b(x)\le 0\}$, for boundness.} The initial set is defined as an intersection of semi-algebraic sets:
\begin{equation}\label{eq:initialset1}
    \mathcal{I}:=\bigcap_{i=1}^l\{x\in\mathbb{R}^n: w_i(x)\ge 0\},
\end{equation}
where $w_1(x),\ldots,w_l(x)$ are all polynomial functions. 
\begin{assumption}\label{ass:local}
    $\mathcal{I}$ is a semi-algebraic set, and $\mathcal{I}$ is bounded on the space $\mathbb{R}^n$.
\end{assumption}
The safe set is defined by
\begin{equation}\label{eq:safeset1}
    \mathcal{S}:=\bigcap_{i=1}^o\{x\in\mathbb{R}^n:a_i^\top (x-c)+1\ge 0\},
\end{equation}
where $a_i\in\mathbb{R}^n$, $c\in\mathbb{R}^n$ is a point in the interior of the safe set. The following theorem proposes a convex condition for $b(x)=(x-c)^\top \Omega^{-1}(x-c)-1$ to be a CBF for $(\eqref{eq:linearsys},\mathcal{I},\mathcal{S})$, with $\mathcal{B}^c=\{x\in\mathbb{R}^n:b(x)\le 0\}$ {a control invariant set.} {By $b(x)$ to be a CBF for $(\eqref{eq:linearsys},\mathcal{I},\mathcal{S})$ we mean that there exists $u(x)$ such that $\frac{\partial b(x)}{\partial x}(Ax+Bu(x))\ge 0$ for all $x\in\partial \mathcal{B}^c$ and $\mathcal{I}\subseteq \mathcal{B}^c\subseteq \mathcal{S}$ (Figure \ref{fig:th2}).} We again use $x_c = x-c$ for notational purpose.

Let $c\in\mathbb{R}^n$ be a {constant} vector such that $\mathrm{rank}([B,Ac])=\mathrm{rank}(B)$ as before, and consider the following optimization program.
\begin{subequations}\label{eq:cqcbflem}
\begin{align}
    \min~&\mathrm{Tr}(\Omega)\\
    \mathrm{subject~to}~&0\prec\Omega=\Omega^\top\in\mathbb{R}^n,\label{eq:cqcbflemeq1}\\
    &0\prec R=R^\top\in\mathbb{R}^{n\times n},Y\in\mathbb{R}^{m\times n},\\
    &\sigma_1(x),\ldots,\sigma_l(x)\in\Sigma[x],\\
    &\Omega A^\top+Y^\top B^\top+A\Omega+BY\preceq 0\label{eq:cqcbflemeqa},\\
    &\begin{bmatrix}
    R&I\\I&\Omega
    \end{bmatrix}\succeq 0\label{eq:cqcbflemeqb},\\
    &-x_c^\top Rx_c+1-\sum_{i=1}^l\sigma_i(x)w_i(x)\in\Sigma[x],\label{eq:cqcbflemeqc}\\
    &1-a_i^\top \Omega a_i\ge 0,i=1,\ldots,o,\label{eq:cqcbflemeqd}
\end{align}
\end{subequations}
where $x_c=x-c$. Similarly to program \eqref{eq:dqcbflem} for global design, program \eqref{eq:cqcbflem} is a convex optimization program, since the cost function is linear, and is subject to semi-definite and linear constraints. In the following theorem, we show how to synthesize a CBF $b(x)$ and a feedback safe controller $u(x)$ by this convex program under Assumption \ref{ass:local}.

\begin{figure}[t]
    \centering
    \includegraphics[scale=0.3]{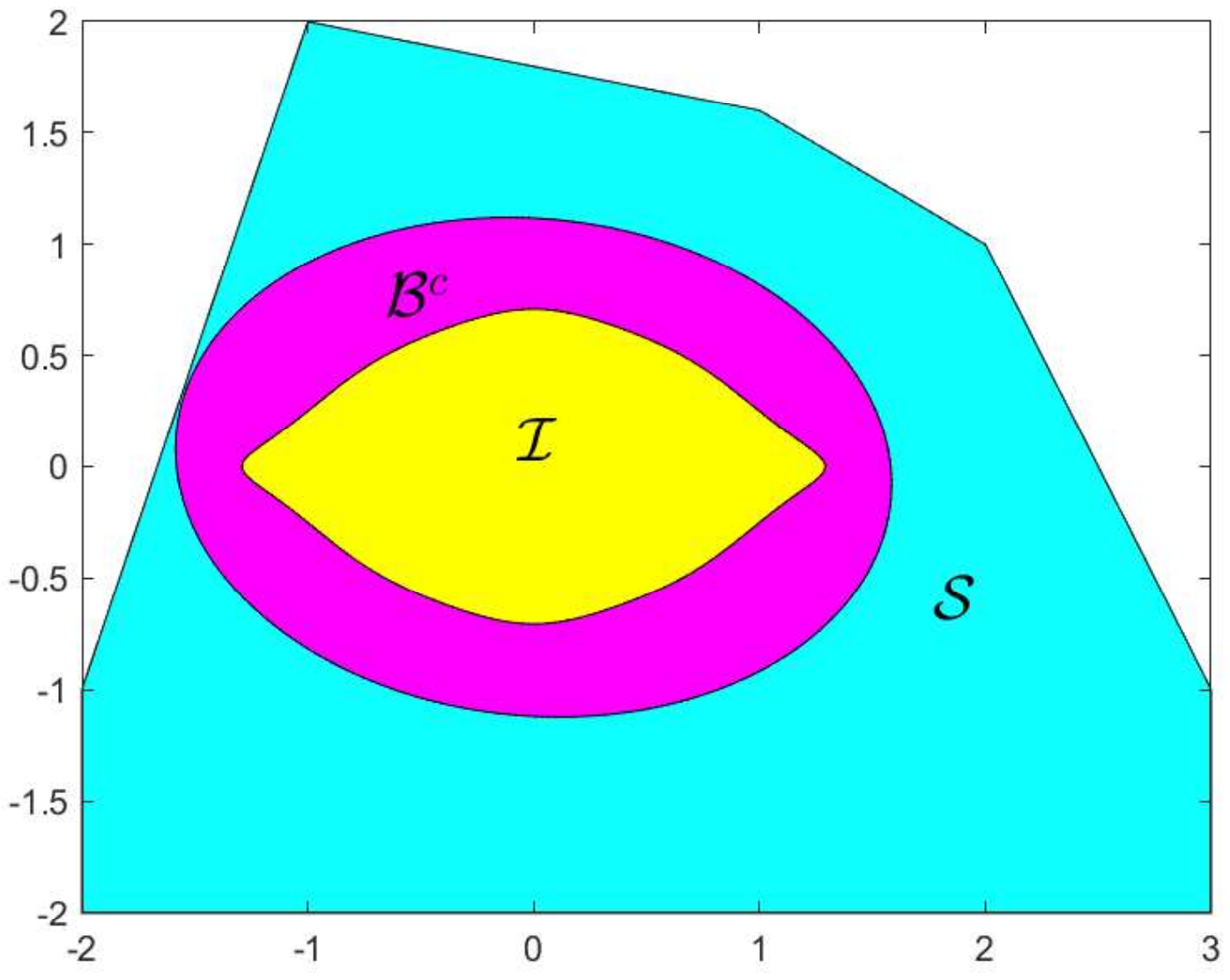}
    \caption{{ Geometric illustration of the the locally constructed CBF $b(x)$ on $\mathbb{R}^2$. The yellow set represents the initial set $\mathcal{I}$, which is bounded on $\mathbb{R}^2$. The blue set represents the safe set $\mathcal{S}$, which is defined by the intersection of half planes on $\mathbb{R}^2$, as in \eqref{eq:safeset1}. The magenta set represents the control invariant set $\mathcal{B}^c:=\{x\in\mathbb{R}^2:b(x)\le 0\}$, which satisfies $\mathcal{I}\subseteq 
    \mathcal{B}^c\subseteq \mathcal{S}$.}}
    \label{fig:th2}
\end{figure}
\begin{theorem}\label{th:cqcbf}
 Consider Assumption \ref{ass:local}, and let $\mathcal{U}=\mathbb{R}^m$. Assume that a solution to \eqref{eq:cqcbflem} exists and is denoted by $\Omega,R,Y,\{\sigma_i(\cdot)\}_{i=1}^l$. Set $u(x)=Y\Omega^{-1}{(x-c)}+d$, where $d\in\mathbb{R}^m$ is such that $Bd+Ac=0$. We then have that
{ 
\begin{enumerate}
    \item $\mathcal{I}\subseteq \mathcal{B}^c\subseteq\mathcal{S}$, where $\mathcal{I}$ is as in \eqref{eq:initialset1}, $\mathcal{B}^c=\{x\in\mathbb{R}^n:b(x)\le 0\}$, $b(x)={(x-c)}^\top \Omega^{-1}{(x-c)}-1$ and $\mathcal{S}$ is as in \eqref{eq:safeset1}.
    \item $\mathcal{B}^c$ is a control invariant set for $\dot x = Ax+Bu(x)$.
\end{enumerate}
}
\end{theorem}
\begin{proof}
The proof that \eqref{eq:cqcbflemeqa} is sufficient for $\mathcal{B}^c$ to be a control invariant set, and {\eqref{eq:cqcbflemeqb} \eqref{eq:cqcbflemeqc}} are sufficient for $\mathcal{I}\subseteq \mathcal{B}^c$ is similar to the proof of Theorem \ref{th:dqcbf}. We only prove that \eqref{eq:cqcbflemeqd} is sufficient and necessary for $\mathcal{B}^c\subseteq \mathcal{S}$.
Using Farkas' lemma \cite{farkas1902theorie}, \cite[Lemma 6.45]{rockafellar2009variational} for affine functions $a_i^\top {(x-c)}+1,i=1,\ldots,o$, and convex quadratic function ${(x-c)}^\top \Omega^{-1}{(x-c)}-1$, we have that $\mathcal{B}^c\subseteq\mathcal{S}$ if and only if for every $i=1,\ldots,o$, there exists $\lambda_i\ge \frac{1}{2}$ such that
\begin{equation*}
    2\lambda_i(a_i^\top x_c +1)-(-x_c^\top \Omega^{-1}x_c+1)\ge 0,\forall x\in\mathbb{R}^n,
\end{equation*}
which is equivalent to
\begin{equation}\label{eq:16}
    \forall i=1,\ldots,o,\exists \lambda_i\ge \frac{1}{2},\mathrm{s.t.}\begin{bmatrix}
        \Omega^{-1}&\lambda_ia_i\\
        \lambda_ia_i^\top & 2\lambda_i-1
    \end{bmatrix}\succeq 0.
\end{equation}
{By Schur complement, \eqref{eq:16} holds if and only if $\Omega\succ 0$, which is true by \eqref{eq:cqcbflemeq1}, and if there exists $\lambda_i\ge \frac{1}{2}$ such that $2\lambda_i-1-\lambda^2a_i^\top \Omega a_i\ge 0$. The discriminant of the quadratic polynomial on the left hand side of the inequality is $4-4a_i^\top \Omega a_i$. Hence, there exists $\lambda_i$ such that $2\lambda_i-1-\lambda_i^2a_i^\top \Omega a_i\ge 0$ if and only if $1-a_i^\top \Omega a_i\ge 0$. Moreover, if $1-a_i^\top \Omega a_i=0$, then $\lambda_i=1\ge \frac{1}{2},$ and if $1-a_i^\top \Omega a_i>0$, then any $\lambda_i\in[1-\sqrt{1-a_i^\top \Omega a_i},1+\sqrt{1-a_i^\top \Omega a_i}]$ satisfies $2\lambda_i-1-\lambda_i^2a_i^\top \Omega a_i\ge 0$. As $1+\sqrt{1-a_i^\top \Omega a_i}>\frac{1}{2},$ we have shown that there exists $\lambda_i\ge\frac{1}{2}$ such that $2\lambda_i-1-\lambda_i^2a_i^\top \Omega a_i\ge 0$ if and only if $1-a_i^\top \Omega a_i\ge 0$, which is \eqref{eq:cqcbflemeqd}. Hence, we conclude the proof.
}
\end{proof}

\subsection{Input Constraints}
In the previous sections for local design we consider the case that $\mathcal{U}=\mathbb{R}^m$. We now extend the local design result to the case that the control authority is limited. Three different types of input constraints are considered: (i) 2-norm bounds, i.e., $\mathcal{U}_1=\{u\in\mathbb{R}^m:||u||_2^2\le \zeta\}$, where $\zeta >0$; (ii) $\infty$-norm bounds, i.e., $\mathcal{U}_2=\{u\in\mathbb{R}^m:||u||_{\infty}\le \sqrt{\zeta}\}$, where $\zeta > 0$; (iii) polytopic bounds, i.e., $\mathcal{U}_3=\{u\in\mathbb{R}^m:Hu\le h\}$, where $H\in\mathbb{R}^{k\times m}$, $h\in\mathbb{R}^k$. 

For $\mathcal{U}=\mathcal{U}_1$, consider the following optimization program {with decision variables $\Omega,R,Y,\sigma_1(x),\ldots,\sigma_l(x),\mu$}:
\begin{subequations}\label{eq:2-norm}
\begin{align}
    \min~\mathrm{Tr}(\overline\Omega)\\
    \mathrm{subject~to}~&\eqref{eq:cqcbflemeq1}-\eqref{eq:cqcbflemeqd},-d^\top d+\zeta-\varepsilon>\mu>0,\\
    &\Pi=\begin{bmatrix}
                \Pi_{11}& \Pi_{12}\\
                \Pi_{12}^\top & I_{n+m+1}
            \end{bmatrix}\succeq 0,\label{eq:2normconstraint}
\end{align}
\end{subequations}
where
\begin{equation*}
        \begin{split}
            \Pi_{11}&=\begin{bmatrix}
                \Omega&Y^\top d&Y^\top\\
                d^\top Y&\mu(-d^\top d+\zeta-\varepsilon)&0\\
                Y&0&\mu
            \end{bmatrix},\\
            \Pi_{12}&=\begin{bmatrix}
                0&0&0\\
                0&\mu&0\\
                0&0&0
            \end{bmatrix},
        \end{split}
        \end{equation*}
and $\varepsilon>0$ is a small constant.
Program \eqref{eq:2-norm} is a convex program which amends program \eqref{eq:cqcbflem} by a new semi-definite constraint \eqref{eq:2normconstraint}.

\begin{lemma}\label{lem:2normconstraint}
    Consider Assumption \ref{ass:local}, and let $\mathcal{U}=\mathcal{U}_1$. Assume that a solution to \eqref{eq:2-norm} exists and is denoted by $\Omega,R,Y,\{\sigma_i(\cdot)\}_{i=1}^l,\mu$. Set $u(x)=Y\Omega^{-1}{(x-c)}+d$, where $d\in\mathbb{R}^n$ is such that $Bd+Ac=0$. We then have that
    { 
    \begin{enumerate}
        \item $\mathcal{I}\subseteq \mathcal{B}^c\subseteq \mathcal{S}$, where $\mathcal{I}$ is as in \eqref{eq:initialset1}, $\mathcal{B}^c=\{x\in\mathbb{R}^n:b(x)\le 0\}$, $b(x)={(x-c)}^\top \Omega^{-1}{(x-c)}-1$ and $\mathcal{S}$ is as in \eqref{eq:safeset1}.
        \item $\mathcal{B}^c$ is a control invariant set for $\dot x = Ax+Bu(x)$ and $u(x)\in\mathcal{U}_1,\forall x \in\mathcal{B}^c$.
    \end{enumerate}
    }
     
\end{lemma} 

\begin{proof}
By Theorem \ref{th:cqcbf}, we have that if $\Omega,R,Y,\{\sigma_i(\cdot)\}_{i=1}^o$ satisfy \eqref{eq:cqcbflemeq1}-\eqref{eq:cqcbflemeqd}, then $\mathcal{B}^c$ is a control invariant set, and $\mathcal{I}\subseteq \mathcal{B}^c\subseteq\mathcal{S}$, and $u(x)$ is a safe controller. We only prove that \eqref{eq:2normconstraint} is sufficient for $u(x)\in\mathcal{U}_1$, for all $x\in\mathcal{B}^c$.
In condition \eqref{eq:2normconstraint}$, \Pi\succeq 0$ is equivalent to
\begin{equation*}
    \begin{split}
        &\Pi_{11}-\Pi_{12}I_{n+m+1}\Pi_{12}^\top\succeq 0\\
        \Longleftrightarrow&
        \begin{bmatrix}
                \Omega&Y^\top d&Y^\top\\
                d^\top Y&\mu(-d^\top d+\zeta-\varepsilon-\mu)&0\\
                Y&0&\mu
            \end{bmatrix}\succeq 0
    \end{split}
\end{equation*}
{ 
By Schur complement, if $\mu>0$, then the last inequality is equivalent to 
\begin{equation}
    \begin{bmatrix}
        \Omega-\mu^{-1}Y^\top Y&Y^\top d\\
        d^\top Y&\mu(-d^\top d+\zeta-\varepsilon-\mu)
    \end{bmatrix}\succeq 0.
\end{equation}
Additionally, $-d^\top d +\zeta -\varepsilon-\mu>0$, then the latter is equivalent to 
\begin{equation}
    \Omega-\mu^{-1}Y^\top Y-\frac{Y^\top dd^\top Y}{\mu(-d^\top d+\zeta-\varepsilon-\mu)}\succeq 0
\end{equation}
or
\begin{equation}
    \mu\Omega-Y^\top Y-\frac{Y^\top dd^\top Y}{-d^\top d+\zeta-\varepsilon-\mu}\succeq 0.
\end{equation}
The matrix remains positive semidefinite if we left- and right- multiply it by $\Omega^{-1}$, thus we obtain (recall that $K=Y\Omega^{-1}$)
}
\begin{align*}
    &\Omega^{-1}\left(\mu \Omega-Y^\top \left(I+\frac{dd^\top}{-d^\top d+\zeta-\varepsilon-\mu}\right)Y\right)\Omega^{-1}\succeq 0,\\
    &\Longleftrightarrow \mu \Omega^{-1}-K^\top K-\frac{K^\top dd^\top K}{-dd^\top +\zeta -\varepsilon-\mu}\succeq 0,\\
    &\Longleftrightarrow
    \begin{bmatrix}
        -K^\top K+\mu\Omega^{-1}&-K^\top d\\
        -d^\top K & -d^\top d+\zeta-\varepsilon-\mu
    \end{bmatrix}\succeq 0,\\
    &\Longleftrightarrow  \\
    &\begin{bmatrix}
        x_c\\1
    \end{bmatrix}^\top \begin{bmatrix}
        -K^\top K+\mu\Omega^{-1}&-K^\top d\\
        -d^\top K & -d^\top d+\zeta-\varepsilon-\mu
    \end{bmatrix}
    \begin{bmatrix}
        x_c\\1
    \end{bmatrix}\ge 0,
\end{align*}
for any $x$.
Writing the product above explicitly, we obtain for any $x\in{\mathbb{R}^n}$:
\begin{align*}
    &-x_c^\top K^\top Kx_c - d^\top Kx_c-x_c^\top K^\top d-d^\top d + \zeta-\varepsilon\\
    +&\mu\left(x_c^\top \Omega^{-1}x_c-1\right) \ge 0,
\end{align*}
which is 
\begin{equation*}
    -\left(u(x)^\top u(x)-\zeta+\varepsilon\right)+\mu({(x-c)}^\top \Omega^{-1}{(x-c)}-1)\ge 0.
\end{equation*}
Hence, for any $x$ such that $b(x)=0$, we have $x_c^\top \Omega^{-1}x_c-1\le0$, then $u(x)^\top u(x)\le \zeta-\varepsilon\le \zeta$. We conclude the proof.
\end{proof}

Condition \eqref{eq:2normconstraint} is an LMI of dimension $2(n+m+1)$. The dimension of the constraints is twice the equivalent condition $\Pi_{11}-\Pi_{12}I_{n+m+1}\Pi_{12}^\top\succeq 0$, which is however not an LMI due to the term $\mu^2$. One tractable convex relaxation while maintaining a relatively lower dimension is
\begin{equation}
    \begin{bmatrix}
        \Omega&Y^\top d&Y^\top\\
        d^\top Y&\frac{(-d^\top d+\zeta-\varepsilon)^2}{2}&0\\
        Y&0&(\frac{-d^\top d+\zeta-\varepsilon}{2})I_m
    \end{bmatrix}\succeq 0.
\end{equation}
Here $\mu$ takes the value of $\frac{-d^\top d+\zeta-\varepsilon}{2}$, which is the maximizer of $\mu(-d^\top d+\zeta-\varepsilon-\mu)$. 

The non-negative tolerance $\varepsilon$ is introduced for robustness. 
\begin{proposition}\label{pro:statedeviation}
Given a CBF $b(x)$, system $\eqref{eq:linearsys}$, and a control admissible set $\mathcal{U}_1$, for any $x$ such that $b(x)=0$, there exists $\delta(x) > 0$, such that for any $x'\in\mathcal{E}(x,\delta(x))$, $u(x')=Y\Omega^{-1}(x'-c)+d\in\mathcal{U}_1$.
\end{proposition}
\begin{proof}
    Given that $u(x)$ is a continuous function, $||u(x)||_2^2$ is also a continuous function. Therefore, for any $x\in\partial \mathcal{B}^c$, there exists $\xi(x)>0$, such that for any $y\in\mathcal{E}(x,\xi(x))$, $||u(y)||_2^2-||u(x)||_2^2\le \frac{\varepsilon}{2}$. From Lemma \ref{lem:2normconstraint} we have that $||u(x)||_2^2\le \zeta-\varepsilon,$ thus $||u(y)||_2^2\le \zeta-\frac{\varepsilon}{2}$. Pick $0<\delta(x)\le \zeta$, we have that for any $x'\in\mathcal{E}(x,\delta(x))$, $||u(x')||\le \zeta-\frac{\varepsilon}{2}$. Hence, $u(x')\in\mathcal{U}_1$, and we conclude the proof.
\end{proof}

We then deal with the case that $\mathcal{U}=\mathcal{U}_2$. Consider the following optimization program {  with decision variables $\Omega,R,Y,\sigma_1(x),\ldots,\sigma_l(x),\mu_1,\ldots,\mu_m$.}

\begin{subequations}\label{eq:inf-norm}
\begin{align}
    \min~\mathrm{Tr}(\overline\Omega)\\
    \mathrm{subject~to}~&\eqref{eq:cqcbflemeq1}-\eqref{eq:cqcbflemeqd},{ -d^\top d+\zeta-\varepsilon>\mu_i>0,}\\
    &\Pi^i=\begin{bmatrix}
                \Pi_{11}^i& \Pi_{12}\\
                \Pi_{12}^\top & I_{n+m+1}
            \end{bmatrix}\succeq 0,i=1,\ldots,m,\label{eq:infnormconstraint}
\end{align}
    
\end{subequations}
where
\begin{equation*}
        \begin{split}
            \Pi_{11}^i&=\begin{bmatrix}
                \Omega&Y^\top O_i^\top d&Y^\top O_i^\top\\
                d^\top O_iY&{ \mu_i}(-d^\top d+\zeta-\varepsilon)&0\\
                O_iY&0&\mu_i
            \end{bmatrix},\\
            \Pi_{12}&=\begin{bmatrix}
                0&0&0\\
                0&{ \mu_i}&0\\
                0&0&0
            \end{bmatrix},
        \end{split}
        \end{equation*}
$O_i\in\mathbb{R}^{m\times m}$ is an all-zero matrix, with the $i$-th diagonal entry is one, $\varepsilon>0$ is a small constant. Program \eqref{eq:inf-norm} is a convex program which amends program \eqref{eq:cqcbflem} by a new semi-definite constraint \eqref{eq:infnormconstraint}.

\begin{lemma}\label{lem:infnormconstraint}
      Consider Assumption \ref{ass:local}, and let $\mathcal{U}=\mathcal{U}_2$. Assume that a solution to \eqref{eq:inf-norm} exists and is denoted by $\Omega,R,Y,\{\sigma_i(\cdot)\}_{i=1}^l,\mu$. Set $u(x)=Y\Omega^{-1}{(x-c)}+d$, where $d\in\mathbb{R}^n$ is such that $Bd+Ac=0$. We then have that
      { 
      \begin{enumerate}
          \item $\mathcal{I}\subseteq \mathcal{B}^c\subseteq \mathcal{S}$, where $\mathcal{I}$ is as in \eqref{eq:initialset1}, $\mathcal{B}^c=\{x\in\mathbb{R}^n:b(x)\le 0\}$, $b(x)={(x-c)}^\top \Omega^{-1}{(x-c)}-1$ and $\mathcal{S}$ is as in \eqref{eq:safeset1}.
          \item $\mathcal{B}^c$ is a control invariant set for $\dot x = Ax+Bu(x)$ and $u(x)\in\mathcal{U}_2$, $\forall x\in\mathcal{B}^c$.
      \end{enumerate}
      }
\end{lemma}
\begin{proof}
    Similarly to the proof of Lemma \ref{lem:2normconstraint}, we prove that \eqref{eq:infnormconstraint} is sufficient for $u(x)\in\mathcal{U}_2,\forall x\in\mathcal{B}^c$. \eqref{eq:infnormconstraint} is equivalent to 
    \begin{align*}
    &-x_c^\top K^\top O_i^\top O_iKx_c - d^\top O_iKx_c-x_c^\top K^\top O_i^\top d-d^\top d + \zeta\\
    +&\mu_i\left(x_c^\top \Omega^{-1}x_c-1\right)-\varepsilon \ge 0,i=1,\ldots,m,
\end{align*}
which is 
\begin{equation*}
\begin{split}
        -\left(u_i(x)^\top u_i(x)-\zeta-\varepsilon\right)+\mu_i({(x-c)}^\top \Omega^{-1}{(x-c)} -1)\ge 0,\\
        i=1,\ldots,m.
\end{split}
\end{equation*}
Then we have $u_i(x)^\top u_i(x)\le \zeta-\varepsilon\le \zeta,i=1,\ldots,m$, for any $x$ such that $b(x)=0$. Therefore, $||u(x)||_\infty\le \sqrt{\zeta}$. We conclude the proof.
\end{proof}

We then deal with the case that $\mathcal{U}=\mathcal{U}_3=\{u\in\mathbb{R}^m:Hu\le h\}$. Consider the following optimization program
\begin{subequations}\label{eq:1-norm}
\begin{align}
    \min~&\mathrm{Tr}(\overline\Omega)\\
    \mathrm{subject~to}~&\eqref{eq:cqcbflemeq1}-\eqref{eq:cqcbflemeqd},\mu>0,\\
    &\begin{bmatrix}
        \Xi _{11}^i&\Xi_{12}\\
        \Xi_{12}^\top&I_{n+1}
    \end{bmatrix}\succeq 0,
    i=1,\ldots,k,\label{eq:onenormconstraint}
\end{align}
\end{subequations}
where
\begin{equation*}
    \Xi_{11}^i=\begin{bmatrix}
        \Omega &Y^\top H_i^\top\\
        H_iY&\mu(-2H_id+2h_i-\varepsilon),
    \end{bmatrix}
    \Xi_{12} = \begin{bmatrix}
        0&0\\0&\mu
    \end{bmatrix},
\end{equation*}
$\varepsilon>0$ is a small constant. Program \eqref{eq:1-norm} is a convex program which amends program \eqref{eq:cqcbflem} by a new semi-definite constraint \eqref{eq:onenormconstraint}. 

\begin{lemma}\label{lem:onenormconstraint}
Consider Assumption \ref{ass:local}, and let $\mathcal{U}=\mathcal{U}_3$. Assume that a solution to \eqref{eq:1-norm} exists and is denoted by $\Omega,R,Y,\{\sigma_i(\cdot)\}_{i=1}^l,\mu$. Set $u(x)=Y\Omega^{-1}{(x-c)}+d$, where $d\in\mathbb{R}^n$ is such that $Bd+Ac=0$. We then have that
{ 
      \begin{enumerate}
          \item $\mathcal{I}\subseteq \mathcal{B}^c\subseteq \mathcal{S}$, where $\mathcal{I}$ is as in \eqref{eq:initialset1}, $\mathcal{B}^c=\{x\in\mathbb{R}^n:b(x)\le 0\}$, $b(x)={(x-c)}^\top \Omega^{-1}{(x-c)}-1$ and $\mathcal{S}$ is as in \eqref{eq:safeset1}.
          \item $\mathcal{B}^c$ is a control invariant set for $\dot x = Ax+Bu(x)$ and $u(x)\in\mathcal{U}_3$, $\forall x\in\mathcal{B}^c$.
      \end{enumerate}
      }
     
\end{lemma}

\begin{proof} 
Similarly to the proof of Lemma \ref{lem:2normconstraint}, and \ref{lem:infnormconstraint}, we prove that \eqref{eq:onenormconstraint} is sufficient for $u(x)\in\mathcal{U}_3,\forall x\in\mathcal{B}$. \eqref{eq:onenormconstraint} is equivalent to
\begin{equation*}
    \Xi_{11}^i-\Xi_{12}I_{n+1}\Xi_{12}^\top\succeq 0,i=1,\ldots,k,
\end{equation*}
which implies
\begin{align*}
    \Omega-\frac{Y^\top H_i^\top YH_i}{\mu(-2H_id+2h_i-\varepsilon-\mu)}\succeq 0,i=1,\ldots,k.
\end{align*}
Therefore, we have
    \begin{equation*}
        \mu\Omega-\frac{Y^\top H_i^\top H_iY}{-2H_id+2h_i-\varepsilon-\mu}\succeq 0.
    \end{equation*}
    {The matrix remains positive semidefinite if we left- and right- multiply it by $\Omega^{-1}$}, thus we obtain
    \begin{align*}
        &\mu\Omega^{-1}-\frac{K^\top H_i^\top H_iK}{-2H_id+2h_i-\varepsilon-\mu}\succeq 0\\
        \Longleftrightarrow&
        \begin{bmatrix}
            \mu\Omega^{-1}&-K^\top H_i^\top\\-H_i K&-2H_id+2h_i-\varepsilon-\mu
        \end{bmatrix}\succeq 0\\
        \Longleftrightarrow&-2(H_i(Kx_c+d)-h_i)+\mu(x_c^\top \Omega^{-1}x_c-1)-\varepsilon\ge 0.
    \end{align*}
Then we have for every $i=1,\ldots,k$, $H_iu_i(x)\le h_i-\varepsilon< h_i$ for any $x$ such that $b(x)={(x-c)}^\top \Omega^{-1}{(x-c)}-1\le0$, $Hu(x)\le h$. We conclude the proof.
\end{proof}

Similar to the design in Lemma \ref{lem:2normconstraint}, $\varepsilon$ is also introduced in Lemma \ref{lem:infnormconstraint} and \ref{lem:onenormconstraint}. {As a consequence, a robustness property is imposed on the synthesized CBF as in Proposition \ref{pro:statedeviation}.}

\section{Simulation Results}
\label{sec:simulation}
In this section we demonstrate the proposed programs on a linear system with a high relative degree. 
All the examples are coded using MATLAB R2022a, SOSTOOLS-4.03 \cite{sostools}, and SeDuMi-1.3.7 \cite{doi:10.1080/10556789908805766}.
\begin{figure}[t]
      \centering
      \includegraphics[scale=0.3]{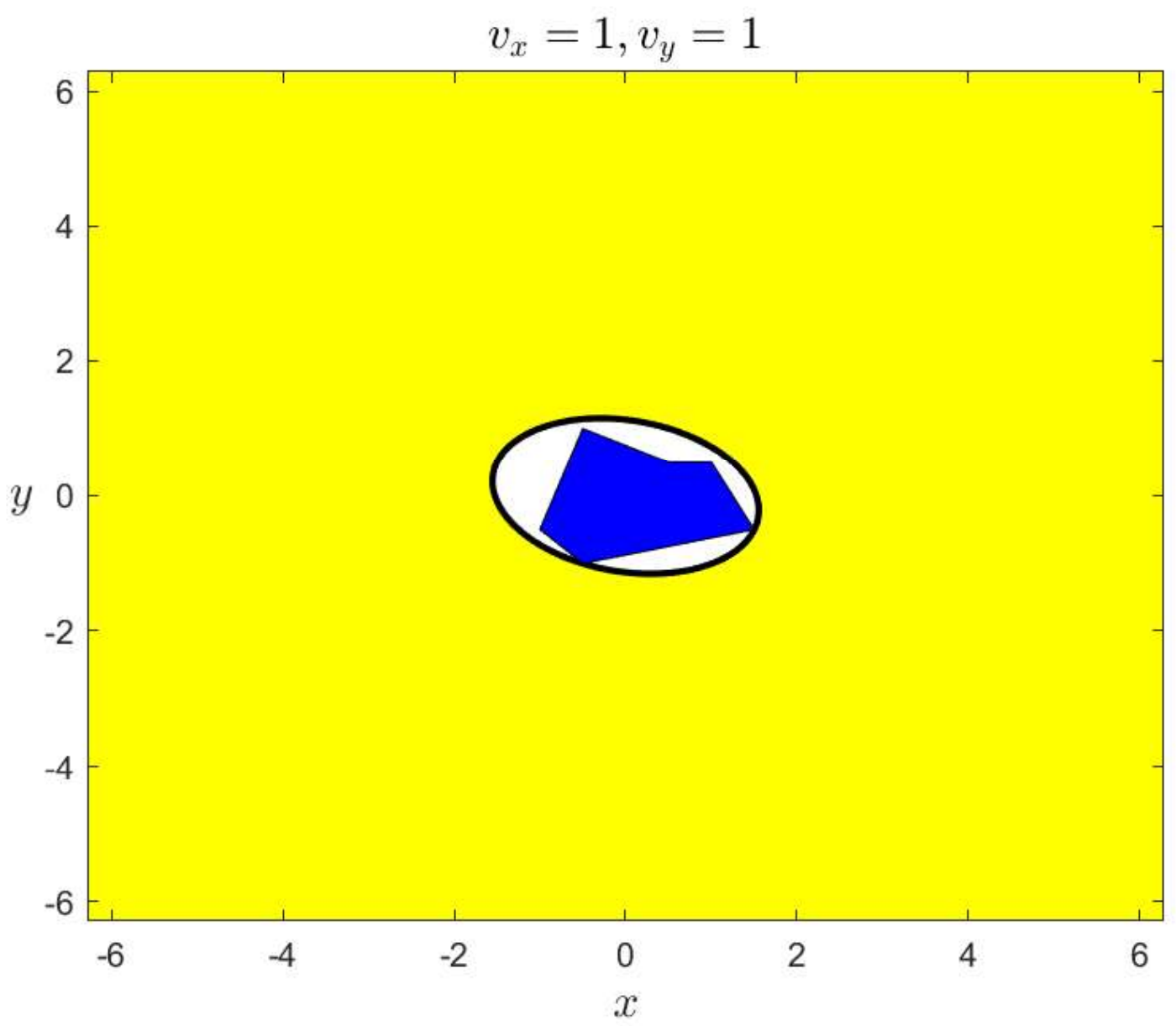}
      \caption{The collision space $\mathcal{S}^c$ is filled by dark blue. Velocity is fixed to be $v_x=1$, $v_y=1$. The control invariant set $\mathcal{B}:=\{x:b(x)\ge 0\}$ is designed by solving the global convex program \eqref{eq:dqcbflem}, and is filled in yellow.}
      \label{fig:global-invariant-omni}
  \end{figure}
In this example, we show how to design CBFs for a linear system with a relative degree. Both the global design and  the local design will be conducted.
Consider an omni-directional vehicle and a collision avoidance problem. The dynamics of the vehicle are

\begin{equation}\label{eq:omnidirectional}
    \begin{bmatrix}
        \dot{x}\\\dot{y}\\\dot{v}_x\\\dot{v}_y
    \end{bmatrix}
    =
    \begin{bmatrix}
        0&0&1&0\\
        0&0&0&1\\
        0&0&0&0\\
        0&0&0&0
    \end{bmatrix}
    {\begin{bmatrix}
        x\\y\\v_x\\v_y
    \end{bmatrix}}
    +
    \begin{bmatrix}
        0&0\\
        0&0\\
        1&0\\
        0&1
    \end{bmatrix}
    \begin{bmatrix}
        a_x\\a_y
    \end{bmatrix},
\end{equation}
where {$[x,y]$} represents the position of the vehicle on the 2-D plane, and {$[v_x,v_y]$} represents the corresponding velocity. The vehicle is controlled by tuning the acceleration denoted by $u={[a_x,a_y]}$ along the two directions. {The position corresponds to $\overline{x}$, while the velocity corresponds to $\underline{x}$ in \eqref{eq:dqcbflem}. }A polytopic obstacle (with five facets) is placed with $\overline c=[0,0]^\top$ be an inside point. Under this configuration, the safe set is a semi-algebraic set, which can be formulated as
\begin{equation*}
    \mathcal{S}:=\bigcup_{i=1}^5\left\{\begin{bmatrix}
        x\\y\\v_x\\v_y
    \end{bmatrix}\in\mathbb{R}^4:a_i^\top \left(\begin{bmatrix}
        x\\y
    \end{bmatrix}-{\overline c}\right)+1\ge 0\right\}.
\end{equation*}
where $a_i\in\mathbb{R}^2$, $i=1,\ldots,5,$ are known vectors.
The collision space $\mathcal{S}^c$ is then a bounded polytope {contains} $c=[0,0,0,0]^\top$. Given that $\mathcal{S}$ is only defined over $[x,y]$, we consider to design a CBF
\begin{equation*}
    b(x)=\begin{bmatrix}
        x\\y
    \end{bmatrix}^\top \overline\Omega^{-1}\begin{bmatrix}
        x\\y
    \end{bmatrix}-\begin{bmatrix}
        v_x\\v_y
    \end{bmatrix}^\top \underline\Omega^{-1}\begin{bmatrix}
        v_x\\v_y
    \end{bmatrix}-1
\end{equation*}
by solving \eqref{eq:dqcbflem}. {We obtain a control barrier function as $b(x)=  2.4104x^2 - 0.67042xy + 1.3229y^2 - 859.4863v_x^2 - 859.4863v_y^2 - 1$ {and a control gain as} $K=\begin{bmatrix}
    369.6&93.6&-0.5&0\\93.6&673.4&0&-0.5
\end{bmatrix}$. We visualize the control invariant set $\mathcal{B}$ and the obstacle $\mathcal{S}^c$ on $\mathbb{R}^2$ by fixing $v_x=1$, $v_y=1$. The result is shown in Figure \ref{fig:global-invariant-omni}.  

\begin{figure}[ht]
     \centering
     \begin{subfigure}{0.4\textwidth}
         \centering
         \includegraphics[width=\textwidth]{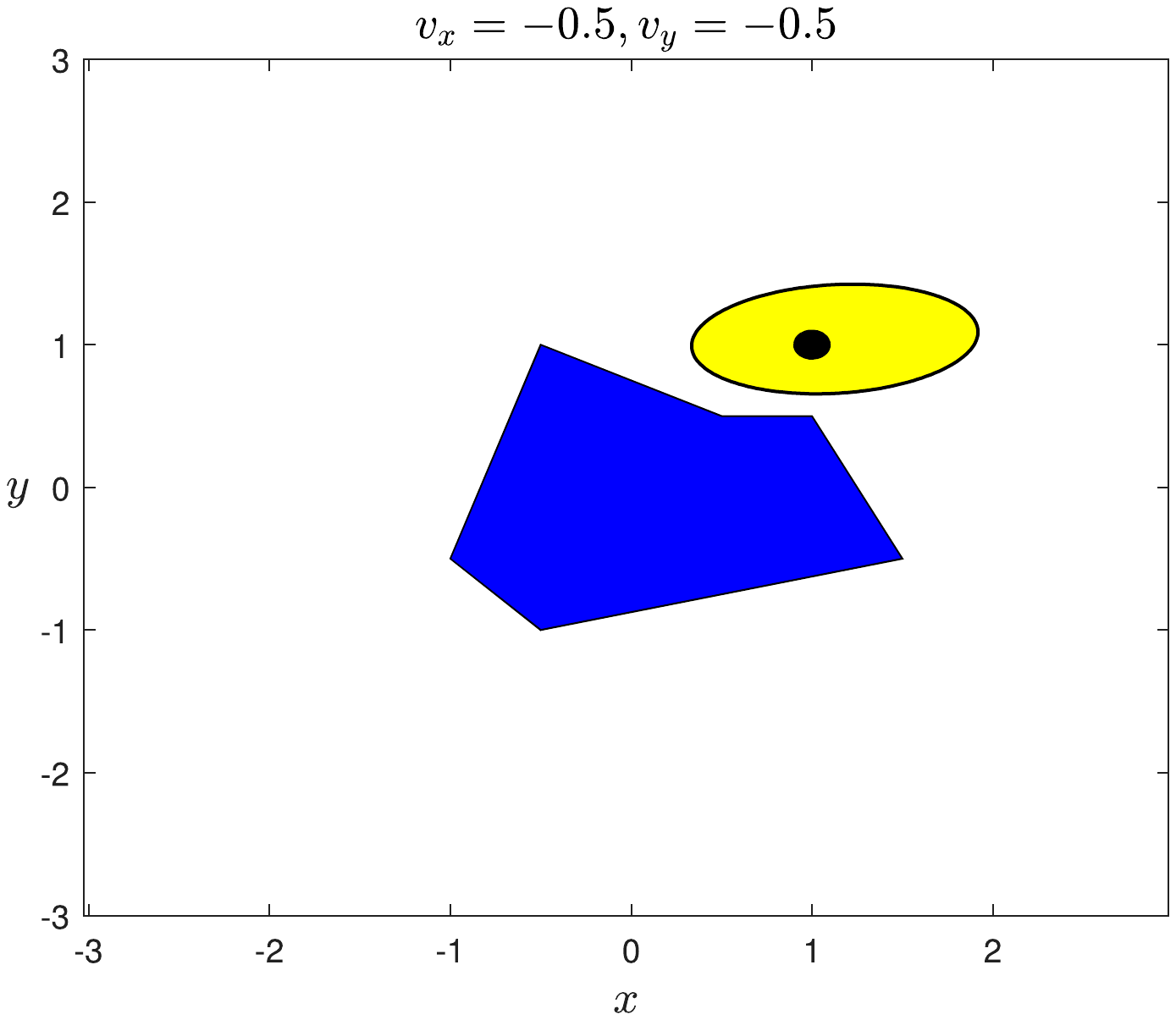}
         \caption{Blue region is the obstacle $\mathcal{S}^c$, yellow region is $\mathcal{B}^c$ and black region is the initial set $\mathcal{I}$.}
         \label{fig:y equals x}
     \end{subfigure}
     \hfill
     \begin{subfigure}{0.4\textwidth}
         \centering
         \includegraphics[width=\textwidth]{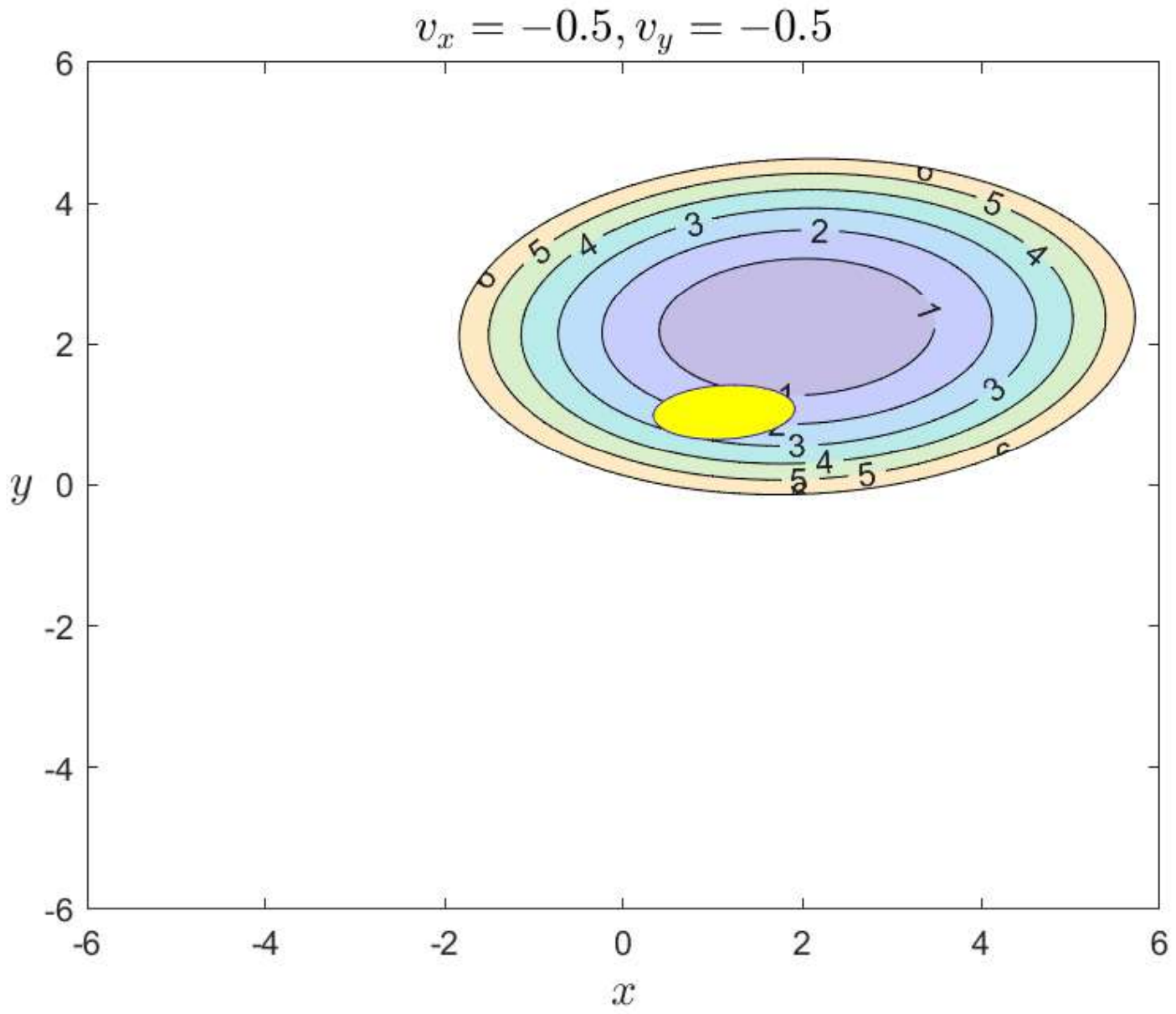}
         \caption{Level sets of $||u(x)||_2^2$. $||u(x)||_2^2\le 4$ for any $x\in\mathcal{B}^c$, thus showing that the input constraint is not violated.}
         \label{fig:three sin x}
     \end{subfigure}
        \caption{Convex invariant set $\mathcal{B}^c$ and state feedback controller $u(x)$ designed by solving the local program \eqref{eq:cqcbflem}. The sets are projected to $\mathbb{R}^2$ by setting $v_x=-0.5,v_y=-0.5$. The yellow region is $\mathcal{B}^c$.}
        \label{fig:omni-local}
\end{figure}

Then we consider a local design. The car is starting from the initial set
\begin{align*}
    \mathcal{I}:=&\left\{\begin{bmatrix}
        x\\y\\v_x\\v_y
    \end{bmatrix}\in\mathbb{R}^4:(x-1)^2+{(y-1)^2}\le 0.01\right.\\
    &\left.\vphantom{\begin{bmatrix}
        x\\y\\v_x\\v_y
    \end{bmatrix}}(v_x+0.5)^2\le 0.1,(v_y+0.5)^2\le 0.1\right\}.
\end{align*}
The acceleration limits are encoded by $a_x^2+a_y^2\le 4.$ By solving the local design program \eqref{eq:cqcbflem}, we obtain $b(x)=0.66903v_x^2 - 0.44567v_xv_y + 0.28291v_xx - 0.80024v_xy + 1.1024
  v_y^2 + 0.23651v_yx + 1.4055v_yy + 1.1198x^2 - 0.58818xy 
  + 4.7544y^2 - 1$ and control gain $K=\begin{bmatrix}
      -1.18&0.41&-0.64&0.21\\
      -0.1&-2.35&-0.09&-1
  \end{bmatrix}.$ The designed control invariant set is $\mathcal{B}^c:=\left\{\begin{bmatrix}
      x\\y\\v_x\\v_y
  \end{bmatrix}\in\mathbb{R}^4:b(x)\le 0\right\}$, and level sets of $||u(x)||_2^2$ are visualized in Figure \ref{fig:omni-local}.
}
\section{Conclusion}\label{sec:conclusion}
In this paper we proposed a method to synthesize a control barrier function and a state feedback controller by solving a single convex program. Our approach considers quadratic control barrier functions and affine state feedback controllers.  
Different types of control input limits can be handled as additional convex constraints to the synthesis program.  We demonstrate the efficacy of our approach on an omni-directional car collision avoidance problem. Future work concentrates towards generalizing the obtained results to allow using higher-relative degree polynomials for the CBF and the controller. We will also consider how to impose input constraint into the global CBF design program using rational polynomial controllers.
\bibliographystyle{ieeetr}
\bibliography{ref.bib}
	\begin{IEEEbiography}[
 {\includegraphics[width=1in,height=1.2in,clip,keepaspectratio]{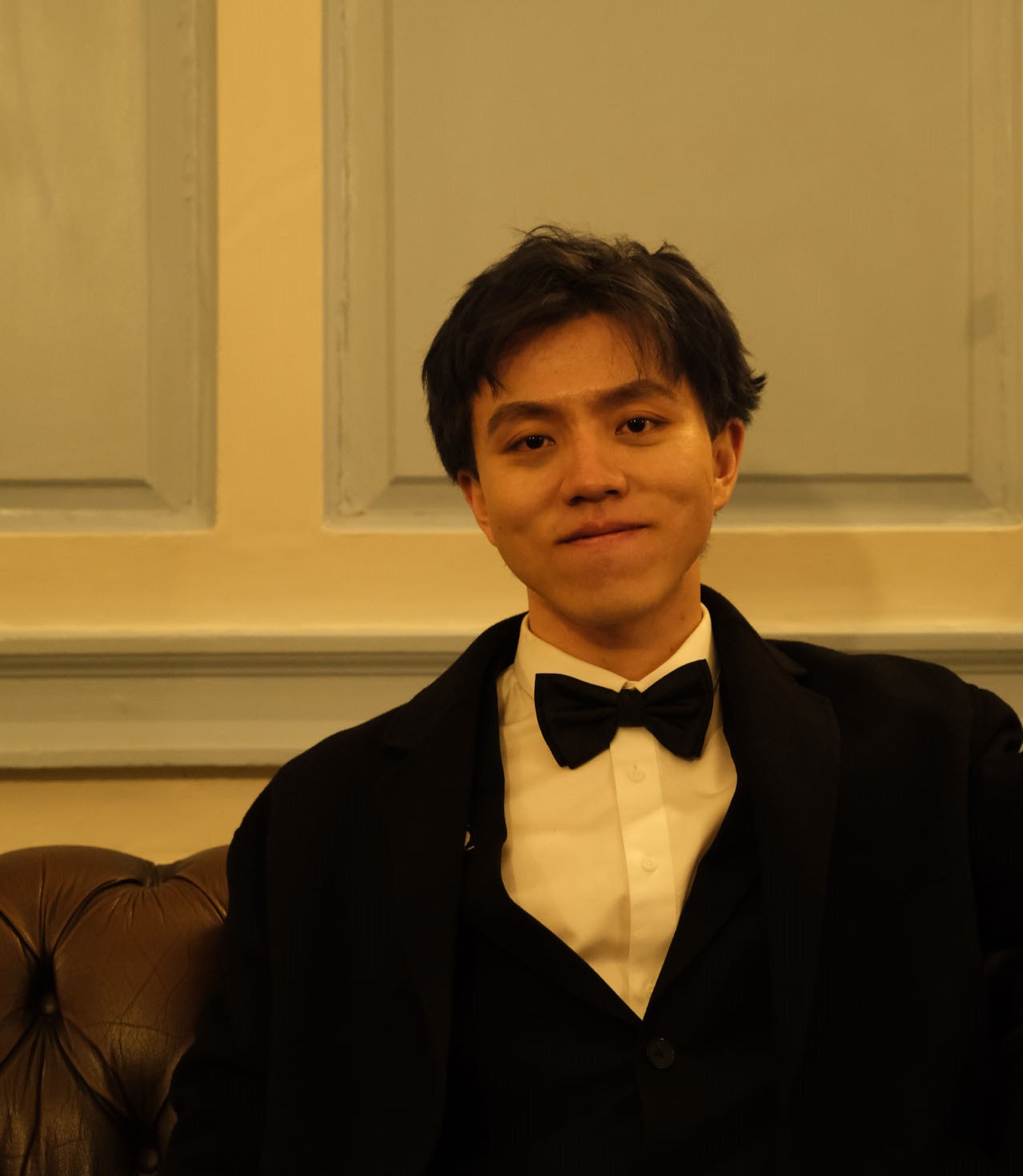}}
		]{Han Wang}
		received B.S. in cyber security from Shanghai Jiao Tong University, China, in 2020. He is currently a Ph.D. student with the Department of Engineering Science at the University of Oxford, Oxford, United Kingdom. His research interests include safe control, data-driven control, and autonomy applications.
	\end{IEEEbiography}
 \vspace{-1cm}
   \begin{IEEEbiography}[
 {\includegraphics[width=1in,height=1.25in,clip,keepaspectratio]{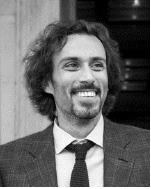}}
		]{Kostas Margellos} received the Diploma in electrical engineering from the University of Patras, Greece, in 2008, and the Ph.D. in control engineering from ETH Zurich, Switzerland, in 2012. He spent 2013, 2014 and 2015 as a postdoctoral researcher at ETH Zurich, UC Berkeley and Politecnico di Milano, respectively. In 2016 he joined the Control Group, Department of Engineering Science, University of Oxford, where he is currently an Associate Professor. He is also a Fellow in AI \& Machine Learning at Reuben College and a Lecturer at Worcester College. 
He is currently serving as Associate Editor in Automatica and in the IEEE Control Systems Letters, and is part of the Conference Editorial Board of the IEEE Control Systems Society and EUCA. His research interests include optimization and control of complex uncertain systems, with applications to energy and transportation networks

	\end{IEEEbiography}
 \vspace{-1cm}
	\begin{IEEEbiography}[
		{\includegraphics[width=1in,height=1.25in,clip,keepaspectratio]{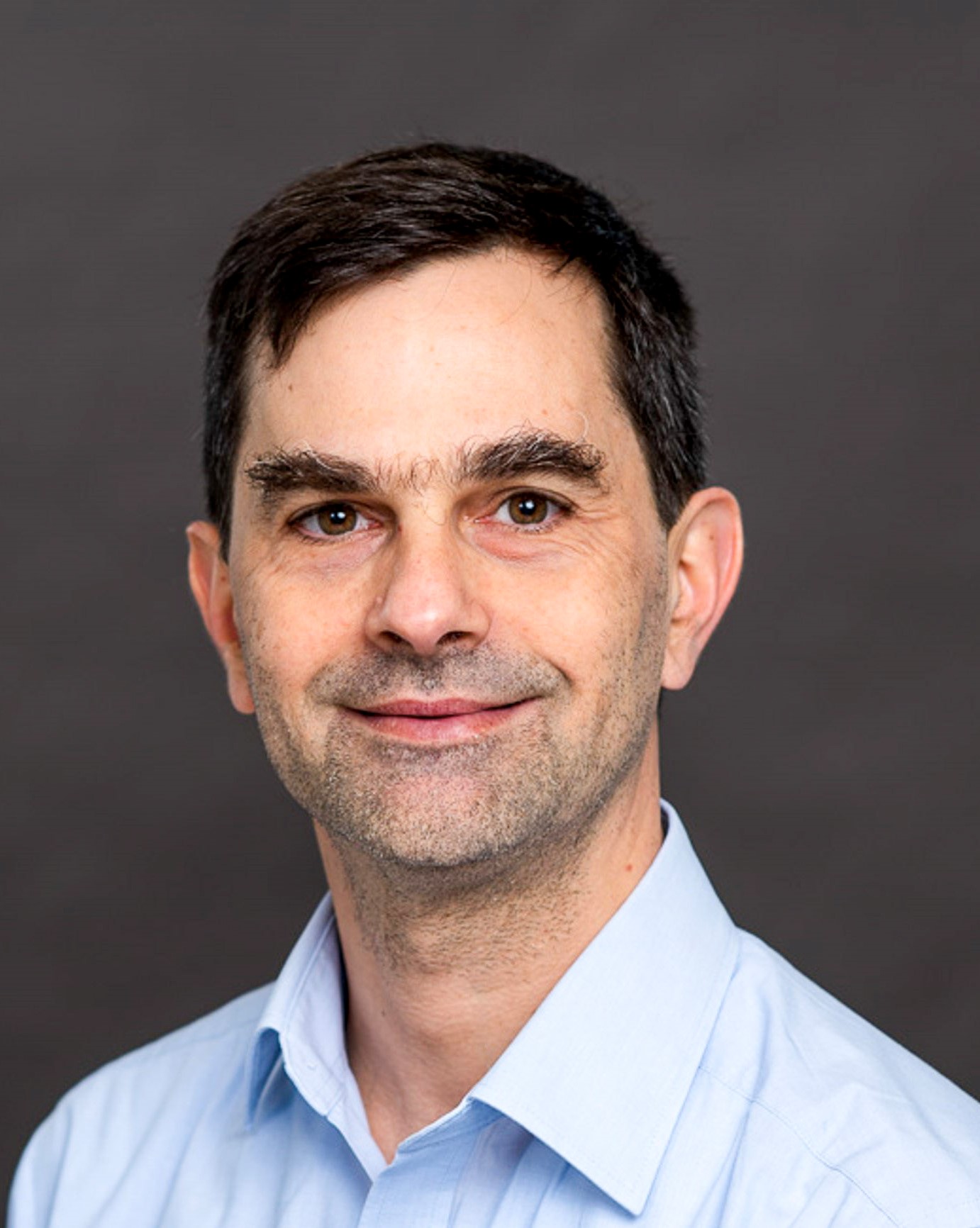}}
		]{Antonis Papachristodoulou} (F'19) received the M.A./M.Eng. degree in electrical and information sciences from the University of Cambridge, Cambridge, U.K., and the Ph.D. degree in control and dynamical systems (with a minor in aeronautics) from the California Institute of Technology, Pasadena, CA, USA. He is currently Professor of Engineering Science at the University of Oxford, Oxford, U.K., and a Tutorial Fellow at Worcester College, Oxford. He was previously an EPSRC Fellow. His research interests include large scale nonlinear systems analysis, sum of squares programming, synthetic and systems biology, networked systems, and flow control.
	\end{IEEEbiography}\vspace{-1cm}
 \begin{IEEEbiography}[
 {\includegraphics[width=1in,height=1.25in,clip,keepaspectratio]{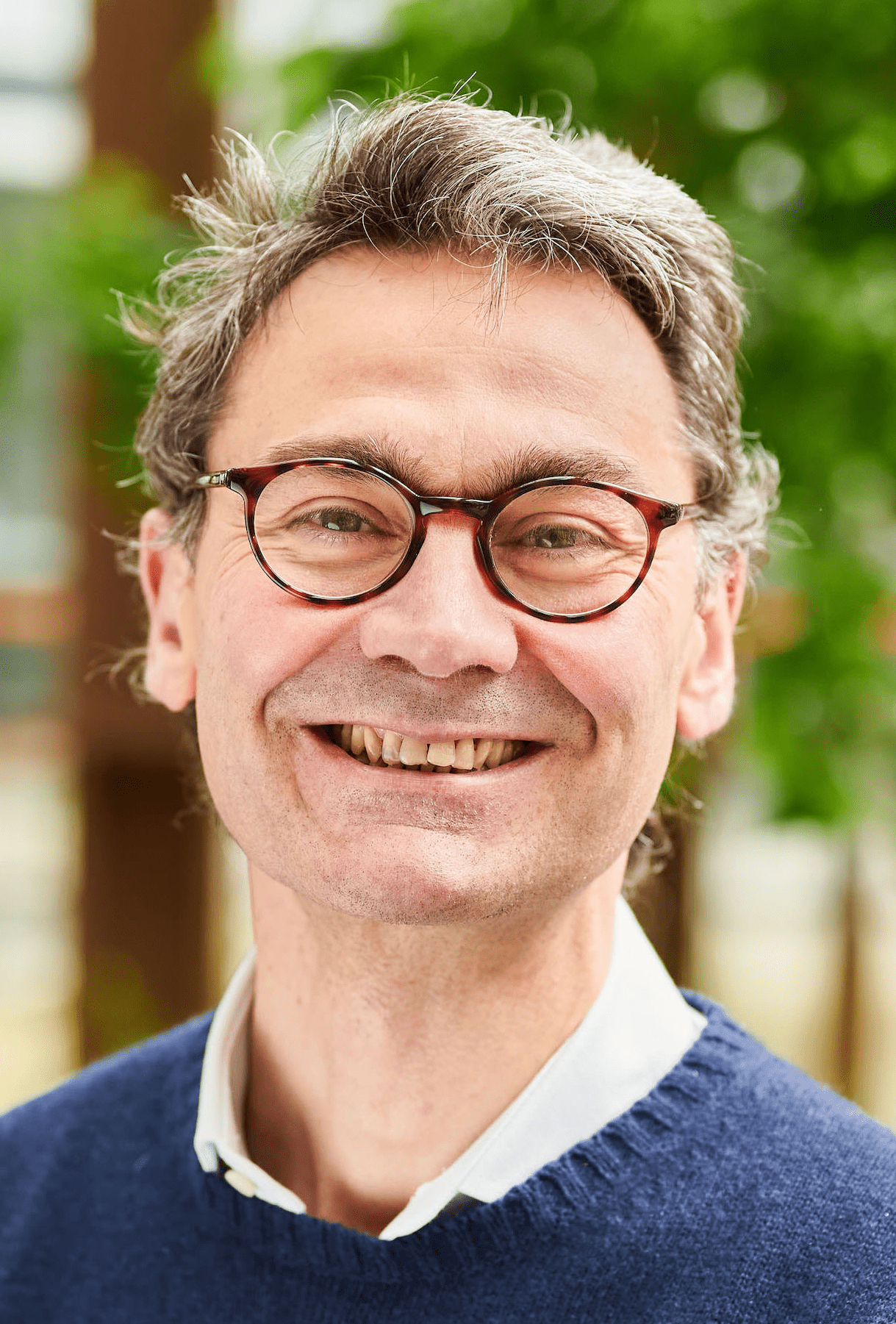}}
		]{Claudio De Persis}  is a Professor with the Engineering and Technology Institute, University of
Groningen, the Netherlands, since 2011. He received
the Laurea and PhD degree in engineering in 1996
and 2000, both from the University of Rome “La
Sapienza”, Italy. He held postdoc positions at Washington University in St. Louis (2000-2001) and Yale
University (2001-2002) and faculty positions at the
University of Rome “La Sapienza” (2002-2009) and
Twente University, the Netherlands (2009-2011). His
main research interest is in automatic control and its
applications.
	\end{IEEEbiography}
 \vspace{-1cm}
\end{document}